%% file: ParticleDynamicsPDECO.tex
\documentclass[final]{siamltex}

\usepackage{amssymb}
\usepackage{amsmath}
\usepackage{mathtools}
\usepackage{latexsym}
\usepackage{graphicx} 
\graphicspath{{Images/}}
\usepackage{url}
\usepackage{xcolor}
\usepackage{multirow}
\usepackage[boxed]{algorithm}
\usepackage{algorithmic} 
\usepackage[warning,debug,nosepfour,autolanguage]{numprint}

\usepackage[framed,numbered,autolinebreaks,useliterate]{mcode}

\newcommand{\adj}{q} 
\newcommand{\hr}{\widehat \rho}

\usepackage{geometry}
\title{Pseudospectral Methods and Iterative Solvers for Optimization Problems from Multiscale Particle Dynamics}

\author{Mildred Aduamoah\thanks{School of Mathematics and Maxwell Institute for Mathematical Sciences, The University of Edinburgh, Edinburgh, EH9 3FD, UK ({\tt maduamoa@ed.ac.uk, b.goddard@ed.ac.uk, j.pearson@ed.ac.uk, J.C.Roden@sms.ed.ac.uk})} 
\and Benjamin D. Goddard\footnotemark[1]
\and John W. Pearson\footnotemark[1]
\and Jonna C. Roden\footnotemark[1]}

\begin{document}
\maketitle

\begin{abstract}
We derive novel algorithms for optimization problems constrained by partial differential equations describing multiscale particle dynamics, including non-local integral terms representing interactions between particles. In particular, we investigate problems where the control acts as an advection `flow' vector or a source term of the partial differential equation, and the constraint is equipped with boundary conditions of Dirichlet or no-flux type. After deriving continuous first-order optimality conditions for such problems, we solve the resulting systems by developing a link with computational methods for statistical mechanics, deriving pseudospectral methods in space and time variables, and utilizing variants of existing fixed-point methods as well as a recently developed Newton--Krylov scheme. Numerical experiments indicate the effectiveness of our approach for a range of problem set-ups, boundary conditions, as well as regularization and model parameters. A key contribution is the provision of software which allows the discretization and solution of a range of optimization problems constrained by differential equations describing particle dynamics.
\end{abstract}



\pagestyle{myheadings}
\thispagestyle{plain}
\markboth{M. ADUAMOAH, B. D. GODDARD, J. W. PEARSON, AND J. C. RODEN}{PSEUDOSPECTRAL METHODS AND SOLVERS FOR OPTIMIZATION OF PARTICLE DYNAMICS}

\section{\label{sec:Intro} Introduction}
In this work we describe a novel approach for tackling partial differential equation (PDE)-constrained optimization problems for systems in which the underlying
dynamics are described by multiscale, interacting particle systems.  Our methods are widely applicable to 
the optimization of many systems described by non-local, non-linear PDEs, some cases
of which have recently received significant attention in the literature~\cite{AcCD10,ABCK,ACFK,AlbiKalise,BonginiButtazo,BKS,CristianiPeri}.

In particular, we aim to provide a link between such optimization problems and state-of-the-art methods in statistical mechanics (known as Dynamic Density
Functional Theory, or DDFT)~\cite{CF05,E79,L10,MT99,RD85,tVLW20,W06,WL07}, before devising numerical methods for such problems using a pseudospectral method in space and time, allowing highly efficient and accurate solution of both the forward and optimization problems~\cite{B01,2DChebClass,NGYSK17,T00}. Having derived first-order optimality conditions using the formal Lagrange method, we modify existing `sweeping', or fixed-point, algorithms~\cite{ACFK,Burger1} to reliably solve such systems, and apply a recently developed Newton--Krylov method~\cite{GP21Software,GP21} to tackle non-linear optimization problems to higher order. The combination of these approaches has not been applied to particle dynamics problems to our knowledge, and enables the accurate resolution of complex optimization problems. We demonstrate how to efficiently implement such methods for problems with both Dirichlet and Robin (no-flux) boundary conditions, provide a number of exact and validation test cases, and accompany the paper with an open-source software implementation~\cite{2DChebClassPDECO}, based on 2DChebClass~\cite{2DChebClass,NGYSK17}.

Our use of pseudospectral methods
has three main advantages over existing implementations: 
(i) due to a novel implementation of spatial convolutions, we are not restricted to periodic domains or the use of Fourier grids and can also tackle convolutions with bounded support; 
(ii) for problems of the types studied here, in particular when the solutions are expected to be smooth and we require accurate solutions, 
pseudospectral methods provide significant computational gains over finite difference or finite element methods;  
(iii) using pseudospectral interpolation in time allows us to move beyond fixed timestepping methods, employing
more accurate and efficient ordinary differential equation (ODE) and differential--algebraic equation (DAE) solvers, and a spectral-in-time Newton--Krylov method.

This paper is structured as follows. In Section~\ref{sec:Background} we provide relevant background to multiscale particle dynamics, pseudospectral
methods, and PDE-constrained optimization.  In Section~\ref{sec:Optimality} we give the associated 
first-order optimality conditions, followed by a description of the numerical methods in Section~\ref{sec:Method}.
The results of our numerical experiments are reported in Section~\ref{sec:Expts}, followed by some concluding 
remarks in Section~\ref{sec:Conc}.

\section{\label{sec:Background} Background}
In this section we detail the necessary background required for the development of our algorithms. In Section \ref{sec:Background_MPS} we describe relevant material on multiscale particle dynamics, in Section \ref{sec:Background_Pseudospectral} we outline pseudospectral methods, in Section \ref{sec:Background_PDECO} we state the PDE-constrained optimization problems of particle dynamics problems that we will consider, and in Section \ref{sec:Background_MFGOC} we survey related work in the area of mean-field optimal control.

\subsection{\label{sec:Background_MPS} Multiscale particle dynamics}
The dynamics of many systems can be accurately described by interacting particles or agents.  Examples 
range in scale from electrons in atoms and molecules~\cite{SO12}, 
through biological cells in tissues~\cite{ACGL12}, up to planets and stars in galaxies~\cite{BT11}.  
Other individual-based models include
animals undergoing flocking and swarming~\cite{YBEM10}, pedestrians walking~\cite{CPT14}, 
or people who interact and thus change their opinions~\cite{L07}.

In principle, such situations can be modelled by differential equations for the `state' (e.g., position, momentum,
opinion) of each individual.  However, the challenge here is that physical systems typically have huge numbers
of particles (e.g., $\sim 10^{25}$ molecules in a litre of water) and, as such, are beyond the treatment of standard
numerical methods, both in terms of storage and processor time.  For $N$ particles, typical algorithms scale
as $N^2$ or $N^3$, which prevents direct computation for more than, say, $\mathcal{O}(10^4)$
particles.  
It is clear from the vast separation of scales between computationally tractable and physically relevant problems
that this issue cannot be overcome through the sequential improvement of computer hardware.

An additional complication of directly solving the dynamics of such systems, e.g., through Newtonian dynamics,
is the sensitive dependence on initial conditions~\cite{LM16}.  For many physical systems, 
it is unreasonable to assume that one knows the exact initial conditions for each particle.  As such, one is interested
not in a particular realization of the dynamics, but rather in an `average' behaviour, which is typical for
the system.
Both of these challenges suggest that it would be prudent to instead study the dynamics through a statistical
mechanics approach, for which one is interested in the macroscopic quantities, rather than individual realizations~\cite{LL94}.
However, this approach comes with its own challenges and drawbacks.  

The first is that, at least without 
additional simplifying approximations, the resulting equations are no easier to solve than the underlying
particle dynamics.  For example, instead of treating the Langevin stochastic differential equation, which 
formally scales computationally as $N^2$, one may treat the corresponding Fokker--Planck (forward Kolmogorov/Smoluchowski)
equation, which is a partial differential equation in $dN$ dimensions, where $d$ is the number of degrees of
freedom of the one-particle phase space (typically 6 when including momentum, and 3 when only considering the
particle positions).  A standard approach would then be to discretize each degree of freedom, reducing the
PDE to a system of coupled ODEs, which  may then, in principle, be solved numerically.  The issue here lies
with the curse of dimensionality: for $M$ points in each degree of freedom, one requires a total of $M^{dN}$ points.
Taking, for the sake of argument, $M=10$ points and $N=10$ particles in three dimensions, 
then the total number of points required is $10^{30}$, which is far too many for a reasonable computation, and far
too few for an accurate solution.

A common approach to overcome this is to use `coarse-graining', which reduces the dimensionality of the system,
generally at the cost of a loss of accuracy or physical effects, and the introduction of unconstrained approximations~\cite{V08}.  
This links to the second challenge, which concerns the multiscale nature of the problem.  In many systems of interest, physically
crucial effects manifest themselves on scales of the particle size, all the way up to the macroscale.  Examples
include volume exclusion of hard particles~\cite{BC12}, biological cellular alignment~\cite{ANHetal10}, 
and nucleation of clusters and clouds~\cite{L12}.
A standard coarse-graining approach would be to ignore effects such as volume exclusion, and treat the whole
system as a bulk, and hence determine quantities such as average densities and orientations~\cite{LL94}.  Whilst this is viable
in homogeneous systems close to equilibrium, it completely fails to capture heterogeneous systems, symmetry breaking,
and many dynamical effects.

However, an extremely efficient and accurate example of coarse-graining which captures such effects 
is Dynamic Density Functional Theory (DDFT)~\cite{CF05,MT99}.
The crucial observation here is that the full $N$-body information in a system is a functional of the 1-body density,
$\rho(\vec{x},t)$ (i.e., the probability of finding any one particle at a given position at a given time).  This is an extension of 
classical density functional theory (DFT) (see, e.g., the early works~\cite{E79,RD85} 
and later reviews~\cite{L10,W06,WL07}), which considers the equilibrium case, and is linked to the celebrated
quantum version~\cite{HK64}.  The main challenge here is that the proof is non-constructive; it is unknown how to 
map from $\rho$ to the full information in the system.  However, in many practical applications, it is $\rho$ itself
that is the quantity of interest.  Hence it is desirable to derive closed equations of motion for the 1-body density,
which is an object in $\mathbb{R}^d$, irrespective of $N$.

The simplest example is the diffusion equation, which corresponds to Brownian motion, and concerns non-interacting
particles; here the reduction to the 1-body density is trivial.  We are instead concerned with systems in which the 
particles interact, e.g., through electrostatic forces, volume exclusion, or exchange of information.  Typical DDFTs
can be thought of as generalized diffusion equations of the form
\begin{equation}
\partial_t \rho (\vec{x},t) = \nabla \cdot \left( \rho \nabla \frac{\delta \mathcal{F}[\rho]}{\delta \rho} \right) = -\nabla \cdot \vec{j}.
\label{eq:DDFT}
\end{equation}
Here $\mathcal{F}$ is the Helmholtz free energy of the system.  For the non-interacting case, at equilibrium, 
it is
\[
\mathcal{F}_{\rm id}[\rho] = \int \rho(\vec{x})( \log \rho(\vec{x}) -1 )~ {\rm d} \vec{x},
\]	
from which it follows that $\nabla \frac{\delta \mathcal{F}_{\rm id} [\rho]}{\delta \rho} = \frac{\nabla\rho}{\rho}$, resulting in
the diffusion equation.  

For more general systems, the exact free energy is unknown (except in the special case of hard rods in one dimension~\cite{TCM08}).
As such, much effort has been devoted to determine accurate approximations of the free energy for a wide range
of systems, but particular focus is given to hard spheres~\cite{R10} and particles with soft interactions~\cite{HM13}; these cases may be
combined in a perturbative manner~\cite{E92}. Here we will focus on a relatively simple DDFT, which closes the equation for 
$\rho$ by considering that the particles are, on average, uncorrelated.  For particles which interact through an even pairwise potential $V_2$,
in an external potential field $V_1$, the (approximate) free energy is modelled by
\[
\mathcal{F}[\rho] = \int \rho(\vec{x})( \log \rho(\vec{x}) -1 )~ {\rm d} \vec{x} + \int V_1(\vec{x}) \rho(\vec{x})~ {\rm d} \vec{x}
+ \frac{1}{2} \int \int \rho(\vec{x}) \rho(\vec{x}') V_2(|\vec{x} - \vec{x}'|) ~ {\rm d} \vec{x} {\rm d} \vec{x}\hspace{0.1em}'.
\]	
This is known as the \emph{mean-field approximation}, which has been shown to be surprisingly accurate for a range of systems~\cite{ACE17},
and is known to be exact
in the limit of dense systems of particles with soft interactions~\cite{MS82}.  We note that this should be considered
as the first stepping stone on a path to treating PDE-constrained optimal control systems for general DDFTs.  Such systems
are highly challenging, not only due to the non-local, non-linear nature of the PDEs, but also due to the complexity of 
the free energy functionals.  For example, Fundamental Measure Theory, which describes the interactions of
systems of hard particles, requires the computation of weighted densities through convolution integrals, followed by
a further integral of a complicated function of these weighted densities~\cite{R10}.  As such, these challenges are postponed
to future work.

A final challenge we will address here is the implementation of (spatial) boundary conditions.  Most physical systems
are constrained in some way, often in a `box' with impassable walls, such that the number of particles is conserved.  
For DDFTs, the corresponding boundary condition is $\vec{j}\cdot \vec{n}=0$ on the boundary, where $\vec{j}$ is the flux, as in \eqref{eq:DDFT}, and
$\vec{n}$ is the unit normal to the boundary.  Whilst this is a standard Neumann boundary condition, we note that the 
difficulty lies in the form of $\vec{j}$; for interacting problems, $\vec{j}$ is non-local and, as such, so is the corresponding boundary
condition.  This results in an equation which is challenging to solve numerically; see Section~\ref{sec:Method}.

\subsection{\label{sec:Background_Pseudospectral} Pseudospectral methods}
There are a number of standard methods for solving DDFT-like problems.  The two most common are the finite element method (FEM)
and pseudospectral methods.  Here we focus on the latter, but note that the algorithm presented below (see Section~\ref{sec:Method}) is general and may
be easily adapted to other numerical methods.  The main challenge in using FEM for DDFT problems lies in their non-locality.
Heuristically, the principal benefits of FEM are that it (i) produces large, but sparse matrices, leading to systems which may be efficiently
solved, for example through the implementation of standard timestepping schemes and carefully-chosen preconditioners (see e.g., \cite{MNN17,PS13,PSW12,RDW10,SW10,Zule11} for PDE-constrained optimization problems); and
(ii) may be applied to complex domains through standard triangulation/meshing routines.
In contrast, for non-local problems such as DDFT the corresponding matrices are not only large, but also dense.  This prevents
the use of standard numerical schemes and significantly increases the computational cost.

Recently, accurate and efficient pseudospectral methods have been developed to tackle these non-local, non-linear DDFTs~\cite{NGYSK17}.
Some details of the implementation will be discussed in Section~\ref{sec:Method}; here we highlight the benefits and challenges.  As is widely 
known~\cite{B01,T00}, pseudospectral methods are extremely accurate for problems with smooth solutions on `nice' 
domains; here `nice' roughly corresponds to domains which may be mapped to the unit square in a simple (e.g.,\ conformal) manner.
They are more challenging to apply on complex domains (although spectral elements can be seen as  a compromise between FEM
and pseudospectral methods~\cite{B01}), and are also of poor accuracy when the solutions are not smooth (the accuracy is 
order $(1/N)^p$ when the solution has $p$ sufficiently nice derivatives~\cite{T00}, but still at the cost of dense matrices).

Their use to treat DDFT problems stems from three main observations: (i) at least in principle, the diffusion term present in
all DDFTs should lead to smoothing of solutions for sufficiently smooth particle interactions; (ii) the pseudospectral matrices are always
dense and, as such, treating non-local terms does not formally affect the numerical cost; (iii) the implementation of non-local 
boundary conditions may be treated via standard algebraic--differential equations solvers, thus removing the need for bespoke
treatments of different boundary conditions.

\subsection{\label{sec:Background_PDECO} PDE-constrained optimization}
In this section we introduce the two main PDE-constrained optimization problem structures that we consider within a multiscale particle dynamics setting. A significant additional complication compared to a standard PDE-constrained optimization problem is the addition of an integral, interaction term. In the following, the terms `flow control' and `source control' refer to the application of the control in the PDE constraint either non-linearly, as a vector field within an advection operator, or linearly, as a scalar source term in the PDE. 

\subsubsection{Flow control problem}
We commence with the following problem involving minimizing a cost functional containing a sum of $L^2$-norm terms within the entire space--time interval $\Omega\times(0,T)$, constrained by a non-linear time-dependent advection--diffusion equation with additional non-local integral term. The control is applied non-linearly in the form of a vector `flow' term:
\begin{align}\label{AdvDiff} 
\begin{split}
\ \min_{\rho,\vec{w}}~~\mathcal{J}(\rho, \vec{w}) \coloneqq  \frac{1}{2}\int_0^T\int_{\Omega}&{}(\rho-\widehat{\rho})^2~{\rm d} \vec{x}{\rm d}t+\frac{\beta}{2}\int_0^T\int_{\Omega}\left\|\vec{w}\right\|^2~{\rm d} \vec{x}{\rm d}t \\
\ \text{s.t.}\quad~\mathcal{D}(\rho,\vec{w})-\nabla_{r}\cdot\mathcal{I}(\rho)&{}=f\hspace{1.75em}\quad\text{on }\Omega\times(0,T), \\
\ \rho&{}=\rho_{0}(\vec{x})\quad\text{at }t=0,
\end{split}
\end{align}
where
\begin{equation*}
\ \mathcal{D}(\rho,\vec{w})=\partial_{t}\rho-\nabla^{2}\rho+\nabla\cdot(\rho\vec{w})-\nabla\cdot(\rho\nabla{}V_{\text{ext}}),\quad\mathcal{I}(\rho)=\kappa\int_{\Omega}\rho(\vec r)\rho(\vec r\,')\vec{K}(\vec r,\vec r\,')~{\rm d}\vec r\,'.
\end{equation*}
Here, $\Omega\subset\mathbb{R}^{d}$, $d\in\{1,2,3\}$, is some given domain with boundary $\partial\Omega$, and $T$ is a prescribed `final time' up to which the process is modelled. The scalar function $\rho$ and the vector-valued function $\vec{w}$ are the \emph{state} and \emph{control variables}, respectively, $\beta>0$ is a given \emph{regularization parameter}, and $\widehat{\rho}(\vec{x},t)$, $V_{\text{ext}}(\vec{x},t)$, $f(\vec{x},t)$, $\rho_{0}(\vec{x})$ are prescribed functions corresponding to the \emph{desired state}, \emph{external potential}, PDE source term, and initial condition, respectively. We highlight that frequently $f(\vec{x},t)=0$,
which results in conservation of mass; one reason we allow the case $f(\vec{x},t)\neq0$ is to enable us to more readily construct analytic test problems for \eqref{AdvDiff}.
Additionally, the non-local integral term models interactions between individual particles, where $\vec{K}$ denotes some vector function.  
We are particularly interested in the case where $\vec{K}$ is odd, i.e., $\vec{K}(\vec r,\vec r\,') = -\vec{K}(\vec r\,',\vec r\,)$; this is the case when
$\vec{K}(\vec r,\vec r\,') = \nabla_r V_2(\vec r-\vec r\,')$ with $V_2(\vec{x}) = V_2(\|\vec{x}\|)$ an even potential.  However, for now we present the results for a general $\vec{K}$.  For $V_2(\|\vec{x}\|)$ decreasing as $\|\vec{x}\| \to \infty$, the integral term models repulsive (attractive) interactions when $\kappa$ is positive (negative).  Of course, much more general choices of $V_2$ are possible. The parameter $\kappa$ models the particle interaction strength. If $\kappa$ is set to zero, the model reduces to a standard non-linear advection--diffusion equation control problem.

We consider two types of boundary conditions imposed on $\rho$, specifically the Dirichlet condition:
\begin{equation}
\ \label{Dirichlet} \rho=c\quad\text{on }\partial\Omega\times(0,T),
\end{equation}
for a given constant $c\in\mathbb{R}$, and the `no-flux type' condition:
\begin{equation}
\ \label{NoFlux} \mathcal{N}(\rho,\vec{w}) + \mathcal{I}(\rho)\cdot\vec{n}=0\quad\text{on }\partial\Omega\times(0,T).
\end{equation}
Here,
\begin{equation*}
\ \mathcal{N}(\rho,\vec{w}) = \frac{\partial\rho}{\partial{}n}-\rho\vec{w}\cdot\vec{n}+\rho\frac{\partial{}V_{\text{ext}}}{\partial{}n},
\end{equation*}
with $\frac{\partial}{\partial{}n}$ denoting the derivative with respect to the normal $\vec{n}$. The latter is a no-flux boundary condition in the classical sense if $f(\vec{x},t)=0$.

\subsubsection{Source control problem}
We also consider the following problem, with an analogous cost functional to the flow control problem, but now with a scalar function for the control variable, which is applied linearly in the form of a PDE source term. This is again minimized subject to a non-linear time-dependent advection--diffusion equation with an additional integral term:
\begin{align}\label{AdvDiff_Linear}
\begin{split}
\ \min_{\rho,{w}}~~\mathcal{J}(\rho, {w}) = \frac{1}{2}\int_0^T\int_{\Omega}{}&(\rho-\widehat{\rho})^2~{\rm d} \vec{x}{\rm d}t+\frac{\beta}{2}\int_0^T\int_{\Omega}w^2~{\rm d} \vec{x}{\rm d}t \\
\ \text{s.t.}\quad\mathcal{D}_l(\rho,w)-\nabla_{r}\cdot\mathcal{I}(\rho)&{}=f\quad\hspace{1.75em}\text{on }\Omega\times(0,T), \\
\ \rho&{}=\rho_{0}(\vec{x})\quad\text{at }t=0,
\end{split}
\end{align}
where
\begin{equation*}
\ \mathcal{D}_l(\rho,w)=\partial_{t}\rho-\nabla^{2}\rho-\nabla\cdot(\rho\nabla{}V_{\text{ext}})-w.
\end{equation*}
This is posed along with the Dirichlet boundary condition \eqref{Dirichlet}, or the `no-flux type' condition:
\begin{equation}
\ \label{NoFlux_Linear} \mathcal{N}_l(\rho)+\mathcal{I}(\rho)\cdot\vec{n}=0\quad\text{on }\partial\Omega\times(0,T),
\end{equation}
where
\begin{equation*}
\ \mathcal{N}_l(\rho)=\frac{\partial\rho}{\partial{}n}+\rho\frac{\partial{}V_{\text{ext}}}{\partial{}n}.
\end{equation*}

We highlight that this paper is focused on fast and effective numerical methods for solving problems of the form \eqref{AdvDiff} and \eqref{AdvDiff_Linear}, as opposed to theoretical questions such as existence, uniqueness, and regularity. We refer to \cite{AcCD10,ACFK,BonginiButtazo,BKS} for discussion of the first two questions for optimization problems of similar structure, such as those arising from mean-field optimal control. For PDE-constrained optimization problems of the structure examined here, it is typical to demand at least $H^1$ regularity in space for the state variable, with at least $L^2$ for the control variable \cite{Troeltzsch}. We note that if `only' this degree of regularity is to be expected, for example if functions such as $\widehat{\rho}$, $f$, $V_{\text{ext}}$, and $\vec{K}$ themselves have low regularity, the pseudospectral methods examined in this work may not perform substantially better than a finite difference method or a FEM, for instance. However, a key feature of the pseudospectral discretization is that it will exploit whatever regularity does exist within the solution, ensuring far superior convergence compared to alternative methods if the solution has a higher degree of regularity. The particle interaction term within the PDE constraints generally introduces dense matrices under discretization, removing one typical advantage of FEMs or finite difference methods over pseudospectral methods, specifically the presence of sparse matrices for local differential operators. Additionally, due to the potentially poor scaling of optimization methods in the number of grid points, it is highly beneficial to use a method which requires fewer discretized points in space and time, motivating the novel methodology presented in the forthcoming sections. We highlight that the software accompanying this work \cite{2DChebClassPDECO} is designed such that the pseudospectral discretization may readily be replaced with matrices arising from alternative spatial discretizations; we base the work on the pseudospectral method specifically as it is a relatively unexplored class of techniques for PDE-constrained optimization problems, including those arising from particle dynamics, which possesses the significant advantages outlined above.

\subsection{\label{sec:Background_MFGOC} Mean-field optimal control}
Mean-field games were first introduced by Lasry and Lions \cite{LASRY2006619, LASRY2006679, LASRY4,Lasry2007}, and independently by Huang, 
Caines, and Malham\'e \cite{Huang1} under the name Nash certainty equivalence, and have been widely studied since then.  
The main challenge over typical PDE-constrained optimization problems arises from the additional non-linear, non-local interaction term. 
Therefore, standard results in optimal control theory cannot readily be applied, and new approaches have to be developed to address theoretical and 
numerical challenges.

The most commonly studied controls are through the flow, e.g., \cite{ACFK}; interaction term, e.g., \cite{Fornasier_2014no2}; or external
agents, e.g., \cite{burger2016controlling}.  A common assumption is that the particle distribution has compact support
\cite{burger2016controlling, burger2019meanfield, fornasier_lisini_orrieri_savare_2019}, which eliminates the need for boundary conditions.
No-flux boundary conditions, which are a principal focus of our work, have been considered in limited settings~\cite{ACFK, carrillo2018no1}.

The two main avenues of research focus on Vlasov-type PDEs arising from the mean-field limit of Cucker--Smale-like~\cite{CuckerSmale1,CuckerSmale2} models of flocking, and Fokker--Planck equations from the same limit of Langevin dynamics.
For the former, Fornasier et al. provided theoretical results on the convergence of the microscopic sparse optimal control problem to a corresponding 
macroscopic problem, using methods of optimal transport and a $\Gamma$-limit argument, proving existence of optimal controls in the mean-field 
setting, see \cite{fornasier_lisini_orrieri_savare_2019, Fornasier_2014no2,Fornasier_2014}.
Additional work on sparse control strategies can be found in \cite{piccoli2014no1}, as well as in the review paper \cite{Fornasier_20161no1}.
In \cite{burger2019meanfield}, convergence results are proved for systems in which the control is applied through interacting, external agents.
For the Fokker--Planck case, analytical research has focused on the derivation of first-order optimality conditions~\cite{ACFK}, 
existence and regularity of optimal controls~\cite{carrillo2019mean}, and convergence of the microscopic optimal control problem to the mean-field limit~\cite{carrillo2018no1,Pinnau_2017}.

In terms of numerical implementations, Strang splitting schemes~\cite{ChengC.Z1976Tiot,gilbertstrang1} are commonly used, 
in particular for control strategies which employ external agents~\cite{burger2016controlling,burger2019instantaneous,Pinnau_2017}, 
in which the numerical results are used to verify convergence in the mean-field limit.  
In \cite{albi2016selective}, different selective control strategies were considered, and an iterative numerical method was chosen, where the interaction term is 
approximated stochastically.  Other approaches involve combining a Chang--Cooper scheme for the forward equation, finite differences for the adjoint
equation, and Monte-Carlo integration~\cite{ACFK} to solve the PDEs.  The optimization step was performed with a sweeping algorithm, with
updates through the gradient equation, which is similar to the gradient descent method in \cite{Burger1}.  Other related numerical work applies
to porous media Fokker--Planck equations~\cite{carrillo2018no1}, as well as the determination of steady state solutions~\cite{albi2014kinetic,Albi_2014no1}.

As described in Section~\ref{sec:Method}, one of our recommended approaches is an optimization scheme that is inspired by existing sweeping algorithms~\cite{ACFK,Burger1}, but with a novel
coupling to pseudospectral methods used to discretize the space and time domains. This composition of methods offers an efficient and accurate solver for a wide class of problems. To our knowledge, it is the first time that pseudospectral methods have been applied to non-local optimal control problems of this form.

\section{\label{sec:Optimality} First-Order Optimality Conditions for Particle Dynamics Models}
In this section we derive the system of PDEs that we need to solve in order to tackle the models \eqref{AdvDiff} and \eqref{AdvDiff_Linear}. In order to obtain first-order optimality conditions for \eqref{AdvDiff} and \eqref{AdvDiff_Linear}, we apply an \emph{optimize-then-discretize method}, meaning we derive appropriate conditions on the continuous level and then consider suitable discretization strategies. The alternative to this approach is the \emph{discretize-then-optimize} method, however we select the former in order to obtain numerical solutions that better reflect the solutions to the continuous first-order optimality conditions. We highlight that an area of active interest in the PDE-constrained optimization community is to construct discretization schemes such that the two approaches coincide (see \cite{CollisHeinkenschloss} for a fundamental example of a problem for which different results are obtained using the two methods). Below we briefly describe how the first-order optimality conditions are formed using the formal Lagrange method, for both flow control and source control problems with different boundary conditions, and refer to \cite{Troeltzsch}, for instance, for a rigorous justification of how such conditions are formed.

\subsection{\label{sec:Optimality_NonlinearDirichlet} Flow control with Dirichlet boundary condition}
We first consider the advection--diffusion constrained optimization problem \eqref{AdvDiff} with the Dirichlet boundary condition \eqref{Dirichlet}. This leads to the continuous Lagrangian:
\begin{equation}
\ \label{Lagrangian} \mathcal{L}(\rho,\vec{w},\adj_1,\adj_2)=\mathcal{J}(\rho,\vec{w})-\int_0^T\int_{\Omega}\left(\mathcal{D}(\rho,\vec{w})-\nabla_{r}\cdot\mathcal{I}(\rho)-f\right)\adj_1~{\rm d} \vec{x}{\rm d}t-\int_0^T\int_{\partial\Omega}(\rho-c)\adj_2~{\rm d}s{\rm d}t,
\end{equation}
where $\adj_1$ and $\adj_2$ correspond to the portions of the \emph{adjoint variable} $\adj$ arising in the interior of the spatial domain $\Omega$ and its boundary $\partial\Omega$, respectively.

To obtain first-order optimality conditions, we first follow the formal Lagrange method for deriving the \emph{adjoint equation} for time-dependent PDE-constrained optimization, see \cite[Chapter 3]{Troeltzsch} for instance. We obtain that the Fr\'{e}chet derivative of $\mathcal{L}$ in the direction of $\rho$ must satisfy $D_{\rho}\mathcal{L}(\bar{\rho},\bar{w},\adj_1,\adj_2)\rho=0$ for all appropriate functions $\rho$. Integrating the relevant terms of \eqref{Lagrangian} by parts and applying Green's formula, any sufficiently smooth $\rho$ such that $\rho(\vec{x},0)=0$ satisfies
\begingroup
\allowdisplaybreaks
\begin{align}
\ \nonumber 0={}&-\int_0^T\int_{\Omega}\left(\mathcal{D}^*(\adj_1,\bar{w})+\widetilde{\mathcal{I}}^*(\bar{\rho},\adj_1)+\widehat{\rho}-\bar{\rho}\right)\rho~{\rm d} \vec{x}{\rm d}t+\int_{\Omega}\adj_{1}(\vec{x},T)\rho(\vec{x},T)~{\rm d} \vec{x}-\int_0^T\int_{\partial\Omega}\adj_2\rho~{\rm d}s{\rm d}t \\
\ \nonumber &\quad+\int_0^T\int_{\Omega}\big[\nabla\cdot(\adj_{1}\nabla\rho)-\nabla\cdot(\rho\nabla{}\adj_{1})-\nabla\cdot(\rho{}\adj_{1}\bar{w})+\nabla\cdot(\rho{}\adj_{1}\nabla{}V_{\text{ext}}) + \nabla \cdot \widetilde{\mathcal{I}}_{\partial \Omega}(\rho,\bar{\rho},\adj_1)\big]~{\rm d} \vec{x}{\rm d}t\\
\ \label{OptCondrho}
\begin{split}
\ ={}&-\int_0^T\int_{\Omega}\left(\mathcal{D}^*(\adj_1,\bar{w})+\widetilde{\mathcal{I}}^*(\bar{\rho},\adj_1)+\widehat{\rho}-\bar{\rho}\right)\rho~{\rm d} \vec{x}{\rm d}t+\int_{\Omega}\adj_1(\vec{x},T)\rho(\vec{x},T)~{\rm d} \vec{x} \\
\ &\quad+\int_0^T\int_{\partial\Omega}\adj_1\frac{\partial\rho}{\partial{}n}~{\rm d}s{\rm d}t+\int_0^T\int_{\partial\Omega}\left[\left(-\frac{\partial{}\adj_1}{\partial{}n}-\adj_{1}\bar{w}\cdot\vec{n}+\adj_{1}\frac{\partial{}V_{\text{ext}}}{\partial{}n} -\adj_2\right)\rho+ \widetilde{\mathcal{I}}_{\partial \Omega}(\rho, \bar{\rho},\adj_1) \cdot \vec{n}\right]~{\rm d}s{\rm d}t,
\end{split}
\end{align}
\endgroup
where
\begin{align*}
\ \mathcal{D}^*(\adj,\vec{w}):={}&-\partial_{t}\adj-\nabla^{2}\adj-\vec{w}\cdot\nabla{}\adj+\nabla{}V_{\text{ext}}\cdot\nabla{}\adj, \\
\ \widetilde{\mathcal{I}}^*(\rho,\adj):={}&\kappa\left(\int_{\Omega}\rho(\vec r \, ')\vec{K}(\vec r,\vec r\,')~{\rm d}\vec r\,'\right)\cdot\nabla_{\vec r}\adj(\vec r)+\kappa\int_{\Omega}\rho(\vec r\,')\vec{K}(\vec r\,',\vec r)\cdot\nabla_{\vec r\,'}\adj(\vec r\,')~{\rm d}\vec r\,',\\
\widetilde{\mathcal{I}}_{\partial \Omega}(\rho,\bar{\rho},\adj):={}&\kappa \adj(\vec r) \rho(\vec r) \int_\Omega \bar{\rho}(\vec r\,')\vec{K}(\vec r, \vec r\,')~{\rm d}\vec r\,' + \kappa \adj(\vec r)\bar{\rho}(\vec r)  \int_\Omega \rho(\vec r \,')\vec{K}(\vec r, \vec r\,')~{\rm d}\vec r\,'.
\end{align*}
Noting first that \eqref{OptCondrho} must hold for all $\rho\in{}C_0^{\infty}(\Omega\times(0,T))$ (i.e., where $\rho(\vec{x},T)$, $\rho(\vec{x},0)$ vanish on $\Omega$, and $\rho$, $\frac{\partial\rho}{\partial{}n}$ vanish on $\partial\Omega$), and observing that $C_0^{\infty}(\Omega\times(0,T))$ is dense on $L^2(\Omega\times(0,T))$, we obtain the adjoint PDE:
\begin{equation*}
\ \mathcal{D}^*(\adj_1,\vec{w})+\widetilde{\mathcal{I}}^*(\rho,\adj_1)=\rho-\widehat{\rho}\quad\text{on }\Omega\times(0,T).
\end{equation*}

Removing the restriction that $\rho(\vec{x},T)$ vanishes on $\Omega$, and arguing similarly, leads to the adjoint boundary condition $\adj_1(\vec{x},T)=0$. From here, we may similarly remove the condition that $\frac{\partial\rho}{\partial{}n}$ vanishes on $\partial\Omega$ to conclude that $\adj_1=0$ on $\partial\Omega\times(0,T)$. Setting the final integral term in \eqref{OptCondrho} to zero then gives the relation between $\adj_1$ and $\adj_2$. Putting all the pieces together, and relabelling $\adj_1$ as $\adj$, we obtain the complete adjoint problem:
\begin{align}\label{Adjoint}
\begin{split}
\ \mathcal{D}^*(\adj,\vec{w})+\widetilde{\mathcal{I}}^*(\rho,\adj)={}&\rho-\widehat{\rho}\quad\text{on }\Omega\times(0,T), \\
\ \adj={}&0\quad\hspace{1.8em}\text{at }t=T, \\
\ \adj={}&0\quad\hspace{1.8em}\text{on }\partial\Omega\times(0,T).
\end{split}
\end{align}

Searching for the stationary point upon differentiation with respect to each component of $\vec{w}$, using similar working as above, gives:
\begin{align*}
\ 0=D_{w_i}\mathcal{L}(\bar{\rho},\bar{w},\adj_1,\adj_2)w_i={}&\beta\int_0^T\int_{\Omega}\bar{w}_{i}w_i~{\rm d} \vec{x}{\rm d}t-\int_0^T\int_{\Omega}\frac{\partial}{\partial{}x_i}(\bar{\rho}w_i)\adj_1~{\rm d} \vec{x}{\rm d}t \\
\ ={}&\beta\int_0^T\int_{\Omega}\bar{w}_{i}w_i~{\rm d} \vec{x}{\rm d}t+\int_0^T\int_{\Omega}\bar{\rho}\frac{\partial{}\adj_1}{\partial{}x_i}w_i~{\rm d} \vec{x}{\rm d}t-\int_0^T\int_{\Omega}\frac{\partial}{\partial{}x_i}(\bar{\rho}\adj_{1}w_i)~{\rm d} \vec{x}{\rm d}t,
\end{align*}
whereupon considering the derivatives with respect to the all entries of $\vec{w}$, and applying Green's formula, leads to the \emph{gradient equation}:
\begin{equation}
\ \label{Gradient} \beta\vec{w}+\rho\nabla{}\adj=\vec{0}.
\end{equation}

To summarize, the complete first-order optimality system for the problem \eqref{AdvDiff} with the Dirichlet boundary condition $\rho=c$ includes the PDE constraint itself (often referred to as the \emph{state equation}), the adjoint problem \eqref{Adjoint}, and the gradient equation \eqref{Gradient}.

Note that the adjoint terms arising from the particle interactions agrees with the representation of the interaction term in~\cite{ACFK}, where $\kappa\vec{K}(\vec r,\vec r\,') = P(\vec r,\vec r\,')(\vec r\,'-\vec r)$. For the special case when $\vec{K}(\vec r, \vec r \,') = \nabla_r V_2(\|\vec r- \vec r\,'\|)$, we have that $\vec{K}$ is an odd function in the sense that $\vec{K}(\vec r, \vec r\,') = - \vec{K}(\vec r\,',\vec r)$
and
\begin{equation*}
\ \widetilde{\mathcal{I}}^*(\rho,\adj)=
\kappa\int_{\Omega}\rho(\vec r\,') \vec{K}(\vec r,\vec r\,')\cdot  \big[ \nabla_{\vec r}\adj(\vec r) - \nabla_{\vec r\,'}\adj(\vec r\,') \big]~{\rm d}\vec r\,'.
\end{equation*}

\subsection{\label{sec:Optimality_NonlinearNoFlux} Flow control with no-flux type boundary condition}
To provide an illustration of how the same working may be applied to problem \eqref{AdvDiff} with the no-flux boundary condition \eqref{NoFlux}, we briefly consider the Lagrangian given by:
\begin{align*}
\ \mathcal{L}(\rho,\vec{w},\adj_1,\adj_2)=\mathcal{J}(\rho,\vec{w})&{}-\int_0^T\int_{\Omega}\left(\mathcal{D}(\rho,\vec{w})-\nabla_{r}\cdot\mathcal{I}(\rho)-f\right)\adj_1~{\rm d} \vec{x}{\rm d}t \\
\ &-\int_0^T\int_{\partial\Omega}\left(\mathcal{N}(\rho,\vec{w})+\mathcal{I}(\rho)\cdot\vec{n}\right)\adj_2~{\rm d}s{\rm d}t.
\end{align*}
Solving $D_{\rho}\mathcal{L}(\bar{\rho},\bar{w},\adj_1,\adj_2)\rho=0$ for all $\rho$ such that $\rho(\vec{x},0)=0$ gives that:
\begin{align*}
\ 0={}&-\int_0^T\int_{\Omega}\left(\mathcal{D}^*(\adj_1,\bar{w})+\widetilde{\mathcal{I}}^*(\bar{\rho},\adj_1)+\widehat{\rho}-\bar{\rho}\right)\rho~{\rm d} \vec{x}{\rm d}t+\int_{\Omega}\adj_1(\vec{x},T)\rho(\vec{x},T)~{\rm d} \vec{x} \\
\ &\quad+\int_0^T\int_{\partial\Omega}\left[\left((\adj_1-\adj_2)\frac{\partial\rho}{\partial{}n}-\frac{\partial{}\adj_1}{\partial{}n}-(\adj_{1}-\adj_{2})\left(\bar{w}\cdot\vec{n}-\frac{\partial{}V_{\text{ext}}}{\partial{}n} \right)\right)\rho + (\adj_1 - \adj_2) \widetilde{\mathcal{I}}_{\partial \Omega}(\rho, \bar{\rho},\adj_1) \cdot \vec{n}\right]~{\rm d}s{\rm d}t.
\end{align*}
Applying the same reasoning as above then leads to the adjoint problem:
\begin{align*}
\begin{split}
\ \mathcal{D}^*(\adj,\vec{w})+\widetilde{\mathcal{I}}^*(\rho,\adj)={}&\rho-\widehat{\rho}\quad\text{on }\Omega\times(0,T), \\
\ \adj={}&0\quad\hspace{1.8em}\text{at }t=T, \\
\ \frac{\partial{}\adj}{\partial{}n}={}&0\quad\hspace{1.8em}\text{on }\partial\Omega\times(0,T),
\end{split}
\end{align*}
along with the state equation as in \eqref{AdvDiff}, and the gradient equation \eqref{Gradient}.

\subsection{\label{sec:Optimality_LinearDirichlet} Source control with Dirichlet boundary condition}
We next consider the problem \eqref{AdvDiff_Linear} with the Dirichlet boundary condition \eqref{Dirichlet}. This leads to the continuous Lagrangian:
\begin{equation*}
\ \mathcal{L}(\rho,w,\adj_1,\adj_2)=\mathcal{J}(\rho,w)-\int_0^T\int_{\Omega}\left(\mathcal{D}_l(\rho,w)-\nabla_{r}\cdot\mathcal{I}(\rho)-f\right)\adj_1~{\rm d} \vec{x}{\rm d}t-\int_0^T\int_{\partial\Omega}(\rho-c)\adj_2~{\rm d}s{\rm d}t.
\end{equation*}

Solving $D_{\rho}\mathcal{L}(\bar{\rho},\bar{w},\adj_1,\adj_2)\rho=0$ for all $\rho$ such that $\rho(\vec{x},0)=0$ gives that:
\begin{align*}
\begin{split}
\ 0={}&-\int_0^T\int_{\Omega}(-\partial_{t}\adj_{1}-\nabla^{2}\adj_{1}+\nabla{}V_{\text{ext}}\cdot\nabla{}\adj_{1}+\widetilde{\mathcal{I}}^*(\bar{\rho},\adj)+\widehat{\rho}-\bar{\rho})\rho~{\rm d} \vec{x}{\rm d}t \\
\ &\quad+\int_{\Omega}\adj(\vec{x},T)\rho(\vec{x},T)~{\rm d} \vec{x}+\int_0^T\int_{\partial\Omega}\adj_1\frac{\partial\rho}{\partial{}n}~{\rm d}s{\rm d}t,
\end{split}
\end{align*}
with further boundary terms which are eliminated through relating $\adj_1$ and $\adj_2$. This then leads to the adjoint problem:
\begin{align}\label{Adjoint_Linear}
\begin{split}
\ \mathcal{D}_l^*(\adj)+\widetilde{\mathcal{I}}^*(\rho,\adj)={}&\rho-\widehat{\rho}\quad\text{on }\Omega\times(0,T), \\
\ \adj={}&0\quad\hspace{1.8em}\text{at }t=T, \\
\ \adj={}&0\quad\hspace{1.8em}\text{on }\partial\Omega\times(0,T),
\end{split}
\end{align}
where
\begin{equation*}
\ \mathcal{D}_l^*(\adj):=-\partial_{t}\adj-\nabla^{2}\adj+\nabla{}V_{\text{ext}}\cdot\nabla{}\adj.
\end{equation*}

Searching for the stationary point upon differentiation with respect to $w$, using similar working as above, gives:
\begin{equation*}
\ D_w\mathcal{L}(\bar{\rho},\bar{w},\adj_1,\adj_2)w=\beta\int_0^T\int_{\Omega}\bar{w}w~{\rm d} \vec{x}{\rm d}t+\int_0^T\int_{\Omega}\bar{w}\adj_1~{\rm d} \vec{x}{\rm d}t,
\end{equation*}
leading to the \emph{gradient equation}:
\begin{equation}
\ \label{Gradient_Linear} \beta{}w+\adj=0.
\end{equation}

To summarize, the complete first-order optimality system for the problem \eqref{AdvDiff_Linear}, with the Dirichlet boundary condition $\rho=c$, includes the PDE constraint itself, the adjoint problem \eqref{Adjoint_Linear}, and the gradient equation \eqref{Gradient_Linear}.

\subsection{Source control with no-flux type boundary condition}
Applying the same working to problem \eqref{AdvDiff_Linear} with no-flux boundary condition \eqref{NoFlux_Linear}, the Lagrangian is given by:
\begin{align*}
\ \mathcal{L}(\rho,w,\adj_1,\adj_2)=\mathcal{J}(\rho,w)&{}-\int_0^T\int_{\Omega}\left(\mathcal{D}_l(\rho,w)-\nabla_{r}\cdot\mathcal{I}(\rho)-f\right)\adj_1~{\rm d} \vec{x}{\rm d}t \\
\ &-\int_0^T\int_{\partial\Omega}\left(\mathcal{N}_l(\rho)+\mathcal{I}(\rho)\cdot\vec{n}\right)\adj_2~{\rm d}s{\rm d}t.
\end{align*}
Applying the same reasoning as above then leads to the adjoint problem:
\begin{align*}
\ \mathcal{D}_l^*(\adj)+\widetilde{\mathcal{I}}^*(\rho,\adj)={}&\rho-\widehat{\rho}\quad\text{on }\Omega\times(0,T), \\
\ \nonumber \adj={}&0\quad\hspace{1.8em}\text{at }t=T, \\
\ \nonumber \frac{\partial{}\adj}{\partial{}n}={}&0\quad\hspace{1.8em}\text{on }\partial\Omega\times(0,T),
\end{align*}
along with the state equation as in \eqref{AdvDiff_Linear}, and the gradient equation \eqref{Gradient_Linear}.

\section{\label{sec:Method} Numerical Method for the Optimization Model}
In this section we describe the structure of our algorithm for the PDE-constrained optimization models under consideration. After describing a pseudospectral method for the PDE constraints (the forward problem), and the adjoint equations, we outline the optimization solvers to be applied numerically, and detail the measures of accuracy that we will employ in our numerical tests. 
We emphasize that the structure of our algorithm is independent of the choice of solvers in each step, for example,
the pseudospectral method in space may be replaced by finite differences or finite elements for problems with non-smooth solutions.  To highlight
this we will describe two different choices of solver for the optimization stage. Additionally, through
the combination of 2DChebClass~\cite{2DChebClass} and a fixed-point or (spectral-in-time) Newton--Krylov solver, 
one may enforce essentially arbitrary boundary conditions, such as
non-local Robin type, with no additional cost to the user.  This contrasts with traditional `boundary bordering' approaches
~\cite{B01}, for which significant analytical work is often required to derive the correct matrices to impose the boundary conditions
(see below). These properties make the approach highly versatile.

\subsection{\label{sec:Method_PseudospectralForward} Pseudospectral method for the forward problem}
As described in Section~\ref{sec:Background_Pseudospectral}, we solve the forward problem using
Chebyshev pseudospectral methods, in particular implemented in {\scshape {\scshape matlab}} using 2DChebClass~\cite{2DChebClass,NGYSK17}. 
The principal novelties of the method concern the computation of convolution integrals and the implementation
of spatial boundary conditions; the boundary conditions in time will be discussed in the following section. 
This makes the method particularly well-suited to problems on finite, non-periodic domains in which
the interaction term involves a convolution on a region with finite support.  Such applications arise in diverse
fields such as hard-sphere DDFT using Fundamental Measure Theory~\cite{NGYSK17,R10,tVLW20}, 
and opinion dynamics~\cite{L07}.

As described in~\cite{NGYSK17}, the convolution integrals are computed in real space, in contrast to many implementations
in which they are computed via Fourier transforms.  
The principal advantage of Fourier methods is that they are
computationally cheap, requiring only fast Fourier transforms and multiplication of functions.  The main
disadvantage is that for finite, non-periodic domains, one needs to pad the domain, which both increases computational cost for no accuracy gain
and introduces difficulties when applying boundary conditions.
Convolution integrals, including those with bounded support, can be implemented by a single matrix--vector multiplication
in the spatial method, with the matrix
precomputed for all time steps.  Use of the physical domain allows efficient implementation of the boundary conditions.  

As is standard, after discretization, in this case through the use of (mapped) Chebyshev pseudospectral points,
the forward PDE(s) are converted into a system of ODEs.  For example, the diffusion equation becomes
\begin{equation}
\frac{{\rm d}}{{\rm d} t} \boldsymbol{\rho} = D_2 \boldsymbol{\rho}, \qquad \mbox{+ IC and BC},
\label{eq:DiscretizedDiffusion}
\end{equation}
where $\boldsymbol{\rho}$ is a vector of values of the solution at each of the Chebyshev points, and $D_2$ is the
Chebyshev second-order differentiation matrix.
In the interior of the domain, this can be solved using standard time-stepping solvers for ODEs.
The challenge lies in imposing the correct spatial boundary conditions.  One standard approach
is to modify the matrix on the right hand side of \eqref{eq:DiscretizedDiffusion} so that the boundary conditions
are automatically satisfied.  This is known as `boundary-bordering'~\cite{B01}.  For simple boundary conditions, such as
homogenenous Dirichlet or (local) Neumann, such an approach is relatively straightforward.  For example, for homogeneous Dirichlet conditions,
assuming that the initial conditions satisfy the boundary conditions, it is sufficient to set the rows
and columns of $D_2$ that correspond to points on the boundary of the domain to zero.  For homogeneous Neumann, there is a similar approach (see~\cite{T00}), which becomes
more involved with more complex right-hand sides of the PDE.  Another approach is to restrict the computation
to interpolants (solutions) which satisfy the boundary conditions; we do not discuss this here as it is highly
non-trivial for the non-linear, non-local problems that we are interested in.

Here we take a more general approach.  The imposition of spatial boundary conditions can be seen as 
extending the discretized system of ODEs to a system of differential--algebraic equations, where the discretized
PDE is solved on the interior of the domain, and the boundary conditions correspond to algebraic equations.
There are various numerical methods for solving such differential--algebraic equations, see e.g.,~\cite{SEK99} for a Runge--Kutta scheme with algebraic constraints, or~\cite{GP21} for a Newton--Krylov
scheme which allows the inclusion of algebraic constraints alongside the PDE. 
The main advantage here is that the numerical method does not have to be explicitly adapted when 
one changes the boundary conditions; one simply has to specify different algebraic constraints that correspond
to the boundary conditions.  In fact, the 2DChebClass code automatically identifies the boundary of various
geometries, allowing a simple implementation of this approach.\vspace{-1.25em}
\begin{lstlisting}[frame=single]  % Start your code-block

% Construction of RHS of PDE 
function drhodt = rhs(rho)
	
% Definition of flux
flux = -grad*rho + w.*[rho;rho] - kappa*[rho;rho].*(Grad*(Conv*rho));
% Definition of right-hand side of PDE
drhodt = -div*flux;
	
% Application of no-flux boundary conditions
drhodt(bound) = normal*flux;
	
end	
\end{lstlisting}

\subsection{\label{sec:Method_PseudospectralPDECO} Pseudospectral method for the adjoint equation}
For the optimization problem, we have a pair of coupled PDEs: the forward PDE with an initial time condition, and the adjoint equation with a final time condition. Due to the inclusion of Laplacians of opposite sign in the two equations, one must be careful when using a standard time-stepping scheme, since one of the equations will be of backward parabolic form, leading to a (possible) lack of well-posedness and
numerical instability.  For example, the optimality system presented in \eqref{AdvDiff}, \eqref{Adjoint}, and \eqref{Gradient}
results in the adjoint equation being unstable `forward in time'. 
One possible approach
is to apply a backward Euler method for the time derivative in the state equation, with the adjoint operator applied to the adjoint equation, whereupon a huge-scale coupled
system of equations is obtained from matrices arising at each time-step. These may be tackled using a preconditioned iterative method, following
methodology in e.g., \cite{PS13,PSW12,SW10}, but note that the systems considered were sparse whereas our systems are dense.  As above,
as well as boundary conditions in time, there are also boundary conditions in space.  
In contrast, in order to utilize our efficient and accurate forward solver, for our fixed-point approach we reverse time in the adjoint problem, resulting in a set of well-posed equations with initial conditions. For this approach, the forward and adjoint equations are coupled non-locally in time; the adjoint equation requires the value of the state variable at later times, so the two equations cannot be solved simultaneously. By contrast, the Newton--Krylov approach allows us to tackle state and adjoint equations simultaneously.

\subsection{\label{sec:Method_Solver} Optimization solver}
Now that we have presented an accurate and efficient numerical scheme for the solution of a set of PDEs 
which includes the forward and adjoint equations, 
the remaining challenges are to: (i) determine a suitable time discretization for the optimality system;
(ii) choose a suitable optimization scheme. For (i), we again choose a Chebyshev pseudospectral scheme (1D in time), which, assuming that the solutions are smooth in time, leads to exponentially accurate interpolation; it is also the foundation of the spectral-in-time Newton--Krylov
scheme presented in \cite{GP21}.  
For (ii), we note that the choice of optimization solver depends strongly on the nature of the solution, and the amount of information available. 
We consider: (a) a general fixed-point or sweeping method~\cite{ACFK,Burger1}, with an adaptive line search framework to determine a mixing rate~\cite{MaPe04}, to solve the system of equations iteratively, which does not require the analytic computation of the Jacobian, and is
also applicable to problems with box constraints as well as other systems for which the regularity of the solution is not sufficient to be exploited by the spectral-in-time nature of the Newton--Krylov approach; 
(b) a higher-order, more efficient Newton--Krylov scheme, which does require the computation of the Jacobian, and could potentially be more
challenging to adapt to more general problems. 
For the fixed-point method, after applying the pseudospectral discretization, we require the solution of a system of algebraic--differential equations. As in Section \ref{sec:Method_PseudospectralForward}, these can be solved using a standard DAE solver. In this paper, the {\scshape matlab} inbuilt ODE solver \texttt{ode15s} is used.
However, our approach is highly modular and it is straightforward to replace our chosen solvers with any other optimization routine,including space or time discretization of the Newton--Krylov approach. 

In the following, we denote the discretized versions of the variables $\rho$, $\adj$, and $\vec{w}$ by $P$, $Q$, and $W$, respectively. Each of these matrices is of the form $A = [\boldsymbol{a_0}, \boldsymbol{a_1}, ... ,\boldsymbol{a_n}]$, where the vectors $\boldsymbol{a_k}$ represent the solutions at the discretized times $k \in \{0,1,...,n\}$, where $n$ is the number of time steps. In particular, the first column of $P$, denoted by $\boldsymbol{\rho_0}$, corresponds to the initial condition $\rho(\vec{x},0)$. If the spatial domain is one-dimensional, $P$, $Q$, and $W$ are of size $N \times (n+1)$, where $N$ is the number of spatial points. In the two-dimensional case, $P$ and $Q$ are of size $(N_1N_2) \times (n+1)$, where $N_j$ is the number of spatial points in the direction of $x_j$. The discretized control $W$ for linear (source) control problems is also $(N_1N_2) \times (n+1)$  dimensional, while it is $(2N_1N_2) \times (n+1)$ dimensional for non-linear (flow) control problems. 

\subsection{\label{sec:Method_Sweeping}Fixed-point, sweeping method}

We first present a first-order fixed-point method, based on~\cite{ACFK,Burger1}, modified to include
a mixing rate which is standard in density functional theory problems of the type considered here~\cite{R10}.  
We emphasize that this method is included for its simplicity and generality; we have also implemented
a higher order Newton-Krylov method -- see Section~\ref{sec:Method_NK}.
The optimization algorithm is initialized with a guess for the control, $W^{(0)}$. Then, in each iteration, denoted by $i$, the following steps are computed:
\vspace{0.1cm}
\begin{enumerate}
	\item Starting with a guess for the control $W^{(i)}$ as input variable, the corresponding state $P^{(i)}$ is found by solving the state equation.
	\item The adjoint, $Q^{(i)}$, is obtained as the solution of the (reversed in time) adjoint equation, using $W^{(i)}$ and $P^{(i)}$ as inputs. Since $P^{(i)}$ contains the solution for all discretized times $k \in \{0,1,...,n\}$, pseudospectral interpolation circumvents issues resulting from the non-local coupling in time, mentioned in Section \ref{sec:Method_PseudospectralPDECO}. 
	\item The gradient equation is solved for the updated control, $W^{(i)}_g$, using the computed $P^{(i)}$, $Q^{(i)}$.
	\item  The convergence of the optimization scheme is measured by computing the error, $\mathcal{E}$,  between $W^{(i)}$ and $W^{(i)}_{g}$; see Section \ref{sec:Method_Validation}.  If $\mathcal{E}$ is smaller than a set tolerance, the algorithm terminates, otherwise we proceed to Step 5.
	\item We update $W^{(i+1)}$ as a linear combination of the current guess $W^{(i)}$, and the value obtained in step 3, $W^{(i)}_{g}$, employing a mixing rate $\lambda \in [0,1]$:
	\begin{equation}\label{FixedPoint}
	W^{(i+1)} = (1-\lambda)W^{(i)} + \lambda W^{(i)}_{g}.
	\end{equation}
\end{enumerate}
Typical values of $\lambda$, which provide stable convergence in the cases we study here, lie between $0.001$ and $0.1$. Note that, while the solutions $P^{(i)}$ and $Q^{(i)}$ change in each iteration, the initial condition $\boldsymbol{\rho_0}$ and final time condition $\boldsymbol{\adj_n}$ remain unchanged throughout the process; the updates are induced by changing $W^{(i)}$.

It is also feasible to vary the mixing rate $\lambda$, based on viewing $W^{(i)}$ as approximate solutions of a fixed-point problem: an unchanged numerical solution for the control variable $\vec{w}$ at successive iterates indicates that a solution of the PDE-constrained optimization problem has been found. Work in \cite{MaPe04} proposes an adaptive line search framework that can determine $\lambda$ which satisfies an Armijo-type condition and hence converge faster to a fixed-point of a system compared to using a constant mixing rate. Based on this, \cite{Iidu16} uses a potential function $E^{(i)}$ defined as follows based on an iterative scheme:
\begin{equation*}
W^{(i)}(\lambda) := W^{(i)} + \lambda d^{(i)}, \quad D^{(i)}(\lambda) := W^{(i)}(\lambda) - W^{(i)}_g(\lambda), \quad E^{(i)}(\lambda) := \Vert 	D^{(i)}(\lambda) \Vert^2.
\end{equation*}
For the fixed-point scheme \eqref{FixedPoint}, we have that $d^{(i)} := W^{(i)}_g - W^{(i)}$. In this notation, $W^{(i)}(\lambda)$ coincides with that of the subsequent fixed-point iterate, $W^{(i+1)}$. As in \cite{Iidu16}, we consider an adaptive mixing rate which at each iteration satisfies the minimization problem of $E^{(i)}$ over $[0,1]$:
\begin{equation}\label{minarmijowolfe}
\text{Find } \lambda^{(i)} \in [0,1] \quad \text{such that} \quad E^{(i)}(\lambda^{(i)}) = \min_{\lambda \in [0,1]} E^{(i)}(\lambda).
\end{equation}
Since the solution of \eqref{minarmijowolfe} cannot be found easily in general, we may seek an approximate minimum which satisfies the Wolfe-type (or Armijo--Wolfe-type) conditions based on those of~\cite{Wolf69,Wolf71}. Specifically, for some $\delta, \sigma \in (0,1)$ with $\delta < \sigma$, we wish that
\begin{eqnarray}
A^{(i)}(\lambda) &:& ~~ \Vert W^{(i)}(\lambda) - W^{(i)}_g(\lambda)\Vert^2 - \Vert W^{(i)} - W^{(i)}_g \Vert^2 < \delta \lambda \langle W^{(i)} - W^{(i)}_g, d^{(i)}  \rangle, \label{armijo}\\
B^{(i)}(\lambda)&:& ~~ \langle W^{(i)}(\lambda) - W^{(i)}_g(\lambda), d^{(i)}  \rangle > \sigma \langle W^{(i)} - W^{(i)}_g, d^{(i)}  \rangle. \label{curvaturecondition}
\end{eqnarray}
It is possible that with only the Armijo-type condition \eqref{armijo} satisfied, the fixed-point algorithm would not achieve reliable convergence. Hence the condition \eqref{curvaturecondition}, based on the curvature condition discussed in \cite[Section 3.1]{Noce06}, is used additionally to ensure that $\lambda^{(i)}$ is not too small and hence unacceptably short steps are ruled out. Note that for the iterative scheme being applied:
\begin{equation*}
\langle  W^{(i)} - W^{(i)}_g, d^{(i)} \rangle = - \Vert W^{(i)} - W^{(i)}_g \Vert^2, \quad \langle  W^{(i)}(\lambda) - W^{(i)}_g(\lambda), d^{(i)} \rangle = - \langle  W^{(i)}(\lambda) - W^{(i)}_g(\lambda), W^{(i)} - W^{(i)}_g \rangle.
\end{equation*}
Based on discussion in \cite{Lewi13}, we select $\delta = 0.3$ and $\sigma = 0.5$ for our tests. The selection of mixing rate $\lambda$ at each fixed-point iteration is then based on Algorithm \ref{alg:armijowolfe}, whilst also ensuring that $\lambda \geq 0.01$. We set $\lambda_0 = 0.2$. We note that it would also be possible to test the classical (Armijo--Wolfe) conditions using the value of $\mathcal{J}(\rho,\vec{w})$ (see \cite{Wolf69,Wolf71} and \cite[Chapter 3]{Noce06}), but the fixed-point approach described here is computationally cheaper and is found to be effective for our problems.

\begin{algorithm}[H]
	\begin{algorithmic}[1]
		\STATE Set $\alpha = 0$, $\beta = \infty$, $\lambda=\lambda_0$
		\WHILE{$A^{(i)}(\lambda)$ or $B^{(i)}(\lambda)$ does not hold}
		\IF{$A^{(i)}(\lambda)$ does not hold}
		\STATE $\beta = \lambda$
		\ELSIF{$B^{(i)}(\lambda)$ does not hold}
		\STATE $\alpha = \lambda$
		\ELSE
		\STATE $\lambda$ is found
		\ENDIF
		\IF{$\beta < \infty$}
		\STATE $\lambda = \dfrac{1}{2}(\alpha + \beta)$
		\ELSE
		\STATE $ \lambda = 2 \alpha$
		\ENDIF
		\ENDWHILE	
	\end{algorithmic}
	\caption{Armijo--Wolfe Algorithm }
	\label{alg:armijowolfe}
\end{algorithm}

\subsection{\label{sec:Method_NK} Newton--Krylov method}

In addition to the first-order fixed-point method, we also wish to consider a higher-order, Newton-type method, with the aim of achieving satisfactory convergence in many fewer iterations than the fixed-point method. The usual disadvantage of such a method is that one typically needs to solve a number of very large linear systems of equations, unless we design a highly efficient discretization procedure. Further, the linear systems are certainly dense for the particle dynamics problems considered, due to the integral particle interaction terms in the problem.

To circumvent this key difficulty, and exploit the faster convergence achieved by higher-order optimization methods, we employ a recently devised Newton--Krylov method for PDE-constrained optimization problems \cite{GP21} (see also \cite{Kell03,KnKe04} for more general descriptions of such methods), and tailor this to the problem at hand by efficiently describing the PDEs and the associated Jacobian on the discrete level, as well as solving the Newton system efficiently. We highlight that such a method has not previously been applied to PDE-constrained optimization problems which involve integral terms, or problems in which the control variable is applied non-linearly. We now briefly describe how the Newton--Krylov method may be applied to both flow control and source control problems. For both problems the state and adjoint equations may be described in the following general form (see \cite{GP21}), by separating the spatial and temporal derivatives in each case:
\begin{align*}
  \mathbf u'(t) ={}& \mathbf{F}(t,\mathbf u,\mathbf v), \qquad \mathbf u(0) = \mathbf u_0\in \mathbb{R}^N, \\
  \mathbf v'(t) ={}& \mathbf{G}(t,\mathbf u,\mathbf v), \qquad \hspace{-0.25em} \mathbf v(T) = \mathbf 0\in \mathbb{R}^N,
\end{align*}
with the vector-valued functions $\mathbf u,\mathbf v: [0,T] \mapsto \mathbb{R}^N$ denoting the state and adjoint variables $\rho$ and $q$ evaluated at each Chebyshev point in the time variable, and $\mathbf u_0$ corresponding to the initial condition $\rho_0(\vec{x})$. The vector functions $\mathbf{F}$ and $\mathbf{G}$ arise from a method of lines discretization of the state and adjoint PDEs at each time-step, and correspond to the following spatial (derivative and linear) terms:
\begingroup
\allowdisplaybreaks
\begin{align*}
\ &\mathbf{F}(t,\mathbf u,\mathbf v)\leftarrow\left\{\begin{array}{rl}
\nabla^{2}\rho+\frac{1}{\beta}\nabla\cdot(\rho^2\nabla{}\adj)+\nabla\cdot(\rho\nabla{}V_{\text{ext}})+\nabla_{r}\cdot\mathcal{I}(\rho)+f & \text{for flow control}, \\
\nabla^{2}\rho+\nabla\cdot(\rho\nabla{}V_{\text{ext}})+\nabla_{r}\cdot\mathcal{I}(\rho)-\frac{1}{\beta}q+f & \text{for source control}, \\
\end{array}\right. \\
\ &\mathbf{G}(t,\mathbf u,\mathbf v)\leftarrow\left\{\begin{array}{rl}
-\nabla^{2}\adj+\frac{1}{\beta}\rho\,|\nabla{}\adj|^2+\nabla{}V_{\text{ext}}\cdot\nabla{}\adj+\widetilde{\mathcal{I}}^*(\rho,\adj)-\rho+\widehat{\rho} & ~\text{for flow control}, \\
-\nabla^{2}\adj+\nabla{}V_{\text{ext}}\cdot\nabla{}\adj+\widetilde{\mathcal{I}}^*(\rho,\adj)-\rho+\widehat{\rho} & ~\text{for source control}. \\
\end{array}\right.
\end{align*}
\endgroup
Note that the gradient equation \eqref{Gradient} (for the flow control problem) or \eqref{Gradient_Linear} (for the source control problem) has been substituted into the state and adjoint equations where applicable.

Following the working in \cite{GP21}, we may then consider approximations $\widetilde{\mathbf u}_k$, $\widetilde{\mathbf v}_k$ to $\mathbf u$, $\mathbf v$ at the $k$th time-step $t_k$, $k\in\{0,1,...,n\}$, and define Chebyshev interpolants $\widetilde{\mathbf u}(t)$, $\widetilde{\mathbf v}(t)$ based on these approximations. The residual functions:
\begin{equation*}
\ \mathbf r_u(t) := \int_0^t \mathbf F(\tau,\mathbf{\widetilde u}(\tau),\mathbf{\widetilde v}(\tau))~{\rm d} \tau - \mathbf{\widetilde u}(t) + \mathbf{\widetilde u}(0), \quad \mathbf r_v(t) := \int_0^t \mathbf G(\tau,\mathbf{\widetilde u}(\tau),\mathbf{\widetilde v}(\tau))~{\rm d} \tau - \mathbf{\widetilde v}(t) + \mathbf{\widetilde v}(0),
\end{equation*}
can then be approximated at each time-step, along with the exact imposition of initial/final-time conditions, to obtain the expressions:
\begin{align*}
\ [ \mathbf{r}_{u,0},\mathbf{r}_{u,1},\ldots,\mathbf{r}_{u,n} ] ={}& [ \mathbf{F}_0,\mathbf{F}_1,\ldots,\mathbf{F}_n] Q + [\mathbf{\widetilde u}_0 - \mathbf{\widetilde u}_0,\mathbf{\widetilde u}_0 - \mathbf{\widetilde u}_1,\ldots,\mathbf{\widetilde u}_0 - \mathbf{\widetilde u}_n], \\
\ [ \mathbf{r}_{v,0},\mathbf{r}_{v,1},\ldots,\mathbf{r}_{v,n} ] ={}& [ \mathbf{G}_0,\mathbf{G}_1,\ldots,\mathbf{G}_n] Q + [\mathbf{\widetilde v}_0 - \mathbf{\widetilde v}_0,\mathbf{\widetilde v}_0 - \mathbf{\widetilde v}_1,\ldots,\mathbf{\widetilde v}_0 - \mathbf{\widetilde v}_n].
\end{align*}
Here, $\mathbf r_{u,k}$ and $\mathbf r_{v,k}$ approximate $\mathbf r_u(t_k)$ and $\mathbf r_v(t_k)$, $\mathbf{F}_k$ and $\mathbf{G}_k$ denote the functions $\mathbf{F}$ and $\mathbf{G}$ evaluated at time $t_k$, and $Q=[q_{i,j}]_{i,j=1,...,n+1}$ is a $(n+1)\times(n+1)$ collocation matrix arising from cumulative integration. Based on this, we then wish to (approximately) solve $\mathbf R = \mathbf 0$, where the global residual function $\mathbf{R}:\mathbb{R}^{2N(n+1)}\to\mathbb{R}^{2N(n+1)}$ is  given by
\begin{equation*}
\mathbf{R} : \left[\begin{array}{c}
 \mathbf{\widetilde u}_0\\
 \mathbf{\widetilde v}_0\\
 \mathbf{\widetilde u}_1\\
 \mathbf{\widetilde v}_1\\
 \mathbf{\widetilde u}_2\\
 \mathbf{\widetilde v}_2\\
 \vdots\\
 \mathbf{\widetilde u}_n\\
 \mathbf{\widetilde v}_n 
\end{array}\right] \mapsto  
\left[\begin{array}{c}
\mathbf{0}\\
\mathbf{0}\\
\sum_{k=1}^{n+1} q_{k,2} \mathbf{F}_k \\
\sum_{k=1}^{n+1} q_{k,2} \mathbf{G}_k \\
\sum_{k=1}^{n+1} q_{k,3} \mathbf{F}_k \\
\sum_{k=1}^{n+1} q_{k,3} \mathbf{G}_k \\
\vdots\\	
\sum_{k=1}^{n+1} q_{k,n+1} \mathbf{F}_k \\
\sum_{k=1}^{n+1} q_{k,n+1} \mathbf{G}_k
\end{array}\right] 
+
\left[\begin{array}{c}
\mathbf{\widetilde u}_0 - \mathbf{u}_0\\
\mathbf{\widetilde v}_n\\
\mathbf{\widetilde u}_0 - \mathbf{\widetilde u}_1\\
\mathbf{\widetilde v}_0 - \mathbf{\widetilde v}_1\\
\mathbf{\widetilde u}_0 - \mathbf{\widetilde u}_2\\
\mathbf{\widetilde v}_0 - \mathbf{\widetilde v}_2\\
\vdots\\	
\mathbf{\widetilde u}_0 - \mathbf{\widetilde u}_n\\
\mathbf{\widetilde v}_0 - \mathbf{\widetilde v}_{n}
\end{array}\right].
\end{equation*}

Applying Newton iteration for this problem leads to an iterative procedure of the form $\mathbf{x}^{(k+1)} = \mathbf{x}^{(k)} - [\mathbf{J}(\mathbf{x}^{(k)})]^{-1}\mathbf{R}(\mathbf{x}^{(k)})$, with $\mathbf{J}$ denoting the Jacobian matrix of the residual function $\mathbf{R}$. Although Jacobian-free Newton--Krylov methods have been studied \cite{KnKe04}, we elect to form the blocks of the Jacobian matrix explicitly due to the availability of this information for the PDE systems under consideration, in order to achieve rapid convergence of the Newton scheme. This requires us to accurately form the functions $\mathbf F$ and $\mathbf G$, as well as the derivatives of these functions in the directions $\mathbf u$ and $\mathbf v$. In the code below, for the two-dimensional flow control problem with Dirichlet boundary conditions, these quantities are denoted \texttt{JFu}, \texttt{JFv}, \texttt{JGu}, \texttt{JGv}, and our software \cite{2DChebClassPDECO} allows us to compute these quantities to spectral accuracy:\vspace{-1.25em}
\begin{lstlisting}[frame=single]  % Start your code-block

% Definition of state and adjoint PDE operators
K1 = @(t,u,v) L + 2/bet * scalarOperator(dotVectors(grad*u,grad*v)) ...
                + gradVextDotGrad(t) + scalarOperator(LapVext(t)) ...
                + kappa * ( dotVectorOperator(grad*(Conv*u),grad) ...
                            + scalarOperator(L*(Conv*u)) );
K2 = @(t,u,v) -1/bet * scalarOperator(u.^2)*L;
K3 = @(t,u,v) -I + 1/bet * scalarOperator(dotVectors(grad*v,grad*v)) ...
                - kappa * ( Dx1*Conv*scalarOperator(Dx1*v) ...
                            + Dx2*Conv*scalarOperator(Dx2*v) );
K4 = @(t,u,v) L - gradVextDotGrad(t) - kappa * dotVectorOperator(grad*(Conv*u),grad);
f = @(t,u,v) z(t);  g = @(t,u,v) uhat(t);

F = @(t,u,v) K1(t,u,v)*u - K2(t,u,v)*v + f(t,u,v);
G = @(t,u,v) K3(t,u,v)*u - K4(t,u,v)*v + g(t,u,v);

% Specification of Jacobians    
JFu = @(t,u,v) L + 2/bet * scalarOperator(u)*scalarOperator(L*v)  ...
                 + 2/bet * scalarOperator(dotVectors(grad*u,grad*v)) ...
                 + 2/bet * scalarOperator(u) * dotVectorOperator(grad*v,grad) ...
                 + gradVextDotGrad(t) + scalarOperator(LapVext(t)) ...
                 + kappa * ( dotVectorOperator(grad*(Conv*u),grad) ...
                            + scalarOperator(L*(Conv*u)) ...
                            + dotVectorOperator((grad*u),grad)*Conv ...
                            + scalarOperator(u)*L*Conv );     
JFv = @(t,u,v) 1/bet * scalarOperator(u.^2)*L ...
               + 2/bet * scalarOperator(u) * dotVectorOperator(grad*u,grad);    
JGu = @(t,u,v) -I + 1/bet * scalarOperator(dotVectors(grad*v,grad*v)) ...
               + kappa * ( dotVectorOperator(grad*v,grad*Conv) ...
                           - Dx1*Conv*scalarOperator(Dx1*v) ...
                           - Dx2*Conv*scalarOperator(Dx2*v) );
JGv = @(t,u,v) -L + 2/bet * scalarOperator(u) * dotVectorOperator(grad*v,grad) ...
                  + gradVextDotGrad(t) ...
                  + kappa * ( dotVectorOperator(grad*Conv*u,grad) ...
                              - Dx1 * Conv * scalarOperator(u) * Dx1 ...
                              - Dx2 * Conv * scalarOperator(u) * Dx2 );
\end{lstlisting}

In more detail, \texttt{Dx1} and \texttt{Dx2} are matrices applying spatial derivatives in each direction, with \texttt{grad} corresponding to the gradient function, and \texttt{L} the spectral discretization of the Laplacian operator. The function \texttt{Conv} applies a convolution integral with the function $\vec{K}$, \texttt{gradVextDotGrad} applies an operator of the form $\nabla{}V_{\text{ext}}\cdot\nabla$, with \texttt{LapVext} evaluating $\nabla^{2}V_{\text{ext}}$ to spectral accuracy. The function \texttt{scalarOperator} forms a scalar function, with \texttt{dotVectors} taking an inner product of two vectors, and \texttt{dotVectorOperator} similarly taking an inner product of the first argument with the second argument applied to a subsequent term. Finally, \texttt{f} and \texttt{g} describe the source term of the state equation $f$ and the desired state $\widehat{\rho}$ within the PDE operators.

For no-flux type boundary conditions the interior and boundary nodes need to be separated within the code, with the Jacobians defined separately for the boundary conditions. We refer to the open-source software \cite{2DChebClassPDECO} (which also makes use of \cite{GP21Software}) for the full implementation with different boundary conditions, as well as for source control problems. By devising routines to compute all derivatives and integration terms to spectral accuracy for particle dynamics systems, we are able to achieve rapid Newton convergence for a range of problems. Having formed the appropriate terms of the Newton system at each iteration, these are solved inexactly using an inner Krylov method, specifically the Generalized Minimal Residual (GMRES) algorithm \cite{gmres}. Column operations may be applied to $\mathbf J$ so that the leading $2N \times 2N$ block of the Jacobian matrix is invertible, at which point the Kronecker-product based preconditioner described in \cite[Section 2.3]{GP21} may be applied.

We believe there are advantages to both the fixed-point and Newton--Krylov methods we have described above. For a range of problems the higher-order Newton--Krylov method is expected to yield more rapid convergence, due to the inclusion of Jacobian information, and the spectral-in-time representation of the residual along with a pseudospectral discretization in space leads to an efficient solver. By contrast, the fixed-point method does not require the Jacobian matrix (or an approximation to it), and is likely to be applicable to more general problems such as those with additional algebraic constraints, which may have lower regularity and therefore may not be amenable to the spectral-in-time approximation. In Section \ref{sec:Expts} we carry out a number of experiments using both fixed-point and Newton--Krylov methods, to demonstrate and compare their effectiveness.

\subsection{\label{sec:Method_Validation} Measures of accuracy}

All errors in Section \ref{sec:Expts} are calculated as a measure of the difference between a variable of interest, $y$, and a reference value $y_R$, e.g., a previous value of $W^{(i)}$, or an analytic solution to a test problem. The error measure $\mathcal{E}$ is composed of an $L^2$ error in space and an $L^\infty$ error in time. We define absolute and relative $L^2$ spatial errors
\begin{equation*}
\mathcal{E}_{Abs}(t) = \left\| y(\vec{x},t) - y_R(\vec{x},t)\right\|_{L^2(\Omega)},\qquad\mathcal{E}_{Rel}(t) = \frac{\left\| y(\vec{x},t) - y_{R}(\vec{x},t)\right\|_{L^2(\Omega)} }{\left\|y_R(\vec{x},t) \right\|_{L^2(\Omega)} + 10^{-10}},
\end{equation*}
where the small additional term on the denominator prevents division by zero, which are used in the full error measure:
\begin{align*}
\mathcal{E} = \max_{t \in [0,T]}\left[\min\left(\mathcal{E}_{Rel}(t), \mathcal{E}_{Abs}(t)\right)\right].
\end{align*}
The minimum between absolute and relative spatial error is taken to avoid choosing an erroneously large relative error, caused by division of one numerically very small term by another.

We have benchmarked the fixed-point scheme against {\scshape matlab}'s inbuilt \texttt{fsolve} function. The latter uses the trust-region-dogleg algorithm, see~\cite{Powell1}, to solve the optimality system of interest. While it is very robust, it is also much slower than the fixed-point method, which works reliably for the types of problems considered in this paper.

\section{\label{sec:Expts} Numerical Experiments}

The optimal control problems \eqref{AdvDiff} and \eqref{AdvDiff_Linear} require inputs in terms of the desired state $\widehat \rho$, the PDE source term $f$, and the external potential $V_{\text{ext}}$, alongside initial and final time conditions for $\rho$ and $\adj$, respectively. Additionally, an initial guess for the control $\vec{w}$ is needed when using the fixed-point method. These are given  
for each of the examples below.
We also require an interaction kernel, which here we fix as
\begin{equation}
	\vec{K}(\vec r,\vec r\,') = \nabla V_2(\vec r- \vec r\,'), \qquad V_2(\vec{x}) = e^{-\left\|\vec{x}\right\|^2}. \label{V2}
\end{equation}
We note that quality of our results is robust with respect to the precise choice of the interaction kernel; the example
here is chosen for illustration.
Interest lies in how the solution to the optimization problems changes upon varying the interaction strength, $\kappa$.
Here we consider three representative values: $\kappa = 0$ (no interaction),  $\kappa = -1$ (attraction), and  $\kappa = 1$ (repulsion).

As a baseline for the cost, we solve the forward PDE using $\vec{w}=\vec{0}$.  We evaluate the associated cost functional $\mathcal{J}$, 
the value of which is denoted by $\mathcal{J}_{uc}$. We then expect that applying the optimization method lowers the value of the cost functional, which we then aim to minimize by optimizing $\vec{w}$, resulting in a cost $\mathcal{J}_c$. This cost
depends on the value of the regularization parameter $\beta$ and it is expected that the norm of the optimal control applied will increase with decreasing $\beta$. When an initial guess for the control is required, i.e., in the fixed-point method, we take $\vec{w} = \vec{0}$, corresponding to the reference system. The Newton--Krylov solver requires an initial guess for the state and adjoint variables at all times. In the examples below, it suffices to choose the initial and final time conditions for the state and adjoint, respectively, as an initial guess at all time points.

In the following examples, the domain considered is $\Omega \times (0,T) = (-1,1)^d \times (0,1)$; our results are
robust to changes in the domain.  The numbers of spatial (Chebyshev) points are $N_1 = N_2 = 20$ for the two-dimensional examples, and $N_1 = N_2 = N_3 = 20$ for the three dimensional example. The number of time points is $n=11$, unless stated otherwise. The absolute and relative tolerances of the ODE solver ({\scshape matlab}'s \texttt{ode15s}~\cite{SEK99}) for the forward problems are set to $10^{-9}$. The tolerance for the Newton--Krylov and fixed-point solvers are $10^{-16}$ and $10^{-4}$, respectively. The mixing parameter $\lambda$ for each iteration of the fixed-point method is determined using Algorithm~\ref{alg:armijowolfe}, with $\delta = 0.3$ and $\sigma = 0.5$.

\subsection{Two-Dimensional Examples}
We now present four examples, applying no-flux and Dirichlet boundary conditions to both flow and source control problems. The precise
initial condition, external potential, and target chosen in each case are given below.  We recall that the two-body interaction is
given by \eqref{V2}. For our first example (flow control with no-flux boundary conditions, Section~\ref{sec:fcn}),
we show results using both Newton--Krylov and fixed-point solvers.  However, since they produce very similar results, as further
validated in Appendix~\ref{appendix:validation}, for the remaining examples we restrict our results to the more efficient
Newton--Krylov scheme.


\subsubsection{Non-linear (flow) control problem with no-flux boundary conditions}\label{sec:fcn}

\begin{figure}[h]
	\centering
	\includegraphics[scale=0.09]{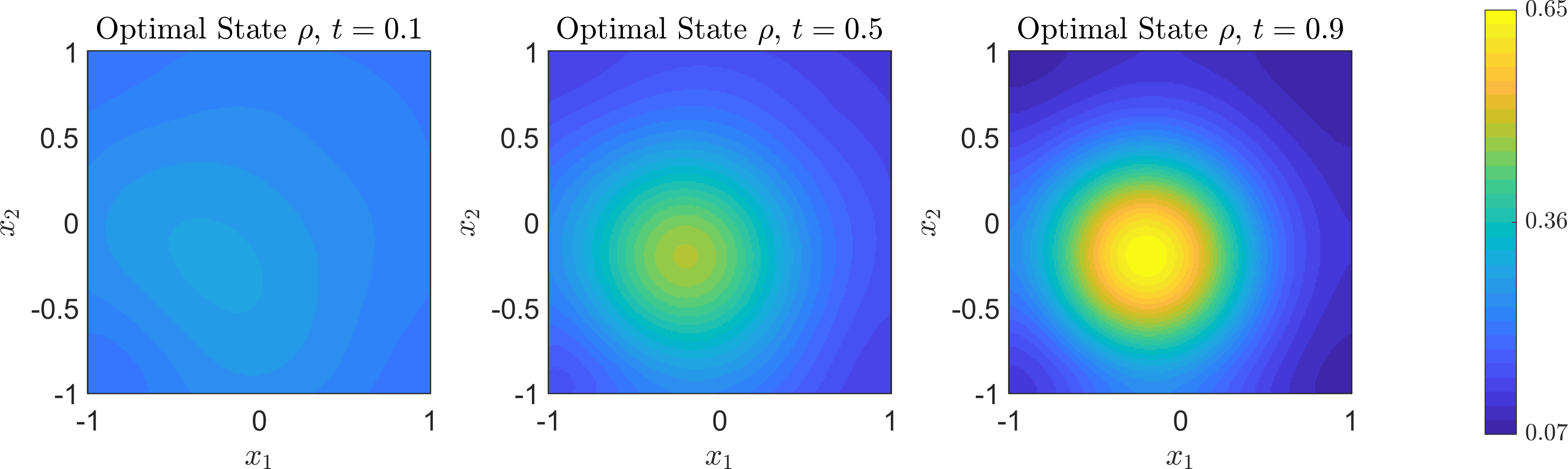}
	\includegraphics[scale=0.09]{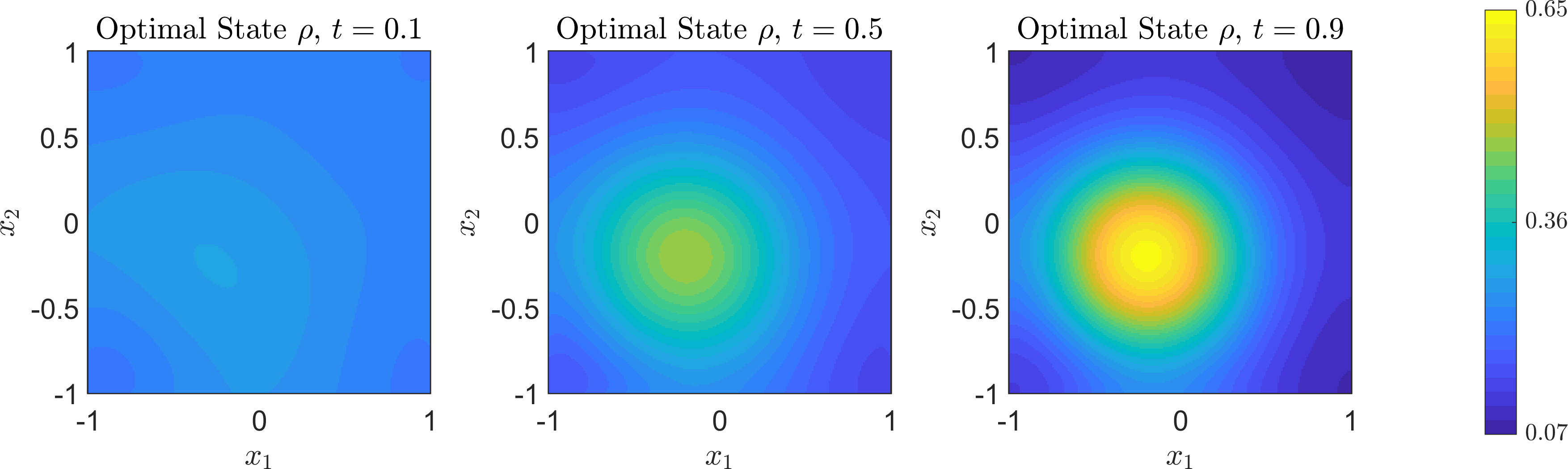}
	\includegraphics[scale=0.09]{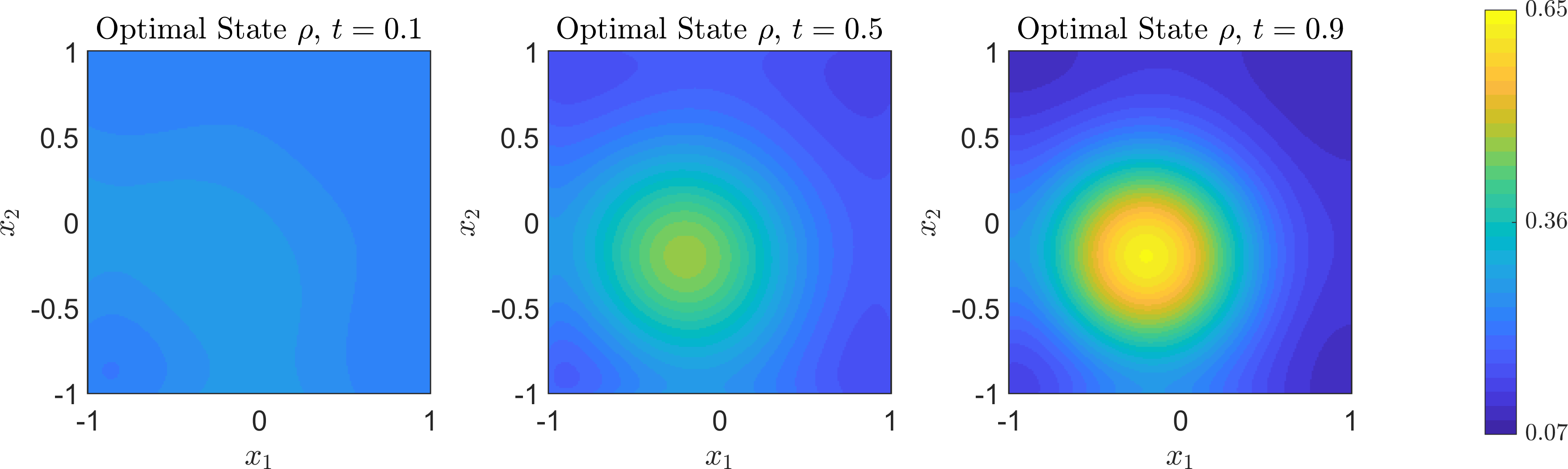}
	\caption{Flow Control, No-Flux: Snapshots of the optimal $\rho$ for different interaction strengths, $\kappa = -1$, $\kappa = 0$, and $\kappa = 1$ (top to bottom), with $\beta = 10^{-3}$.} 
	\label{F3a}
\end{figure}

\begin{figure}[h]
	\centering
	\includegraphics[scale=0.09]{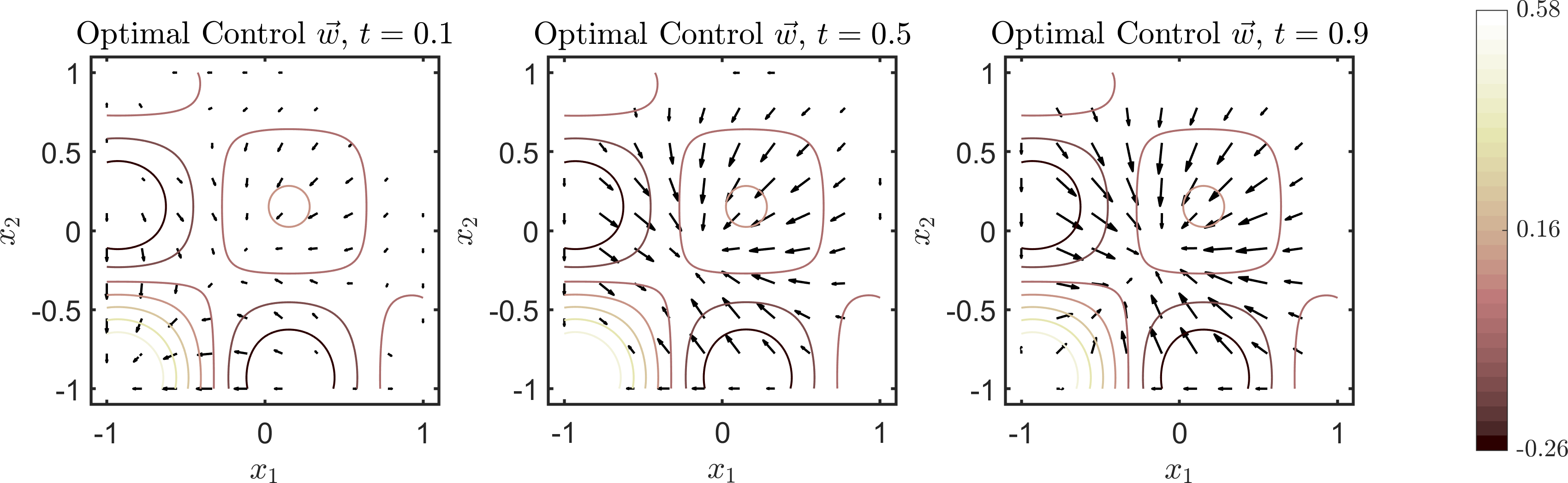}
	\includegraphics[scale=0.09]{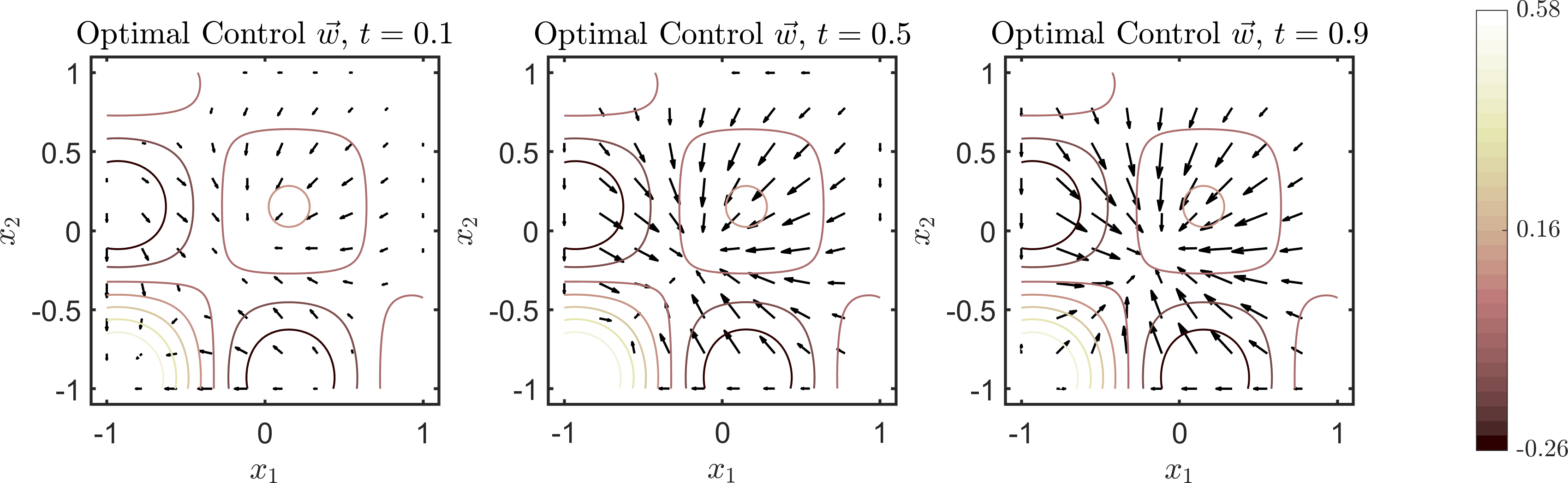}
	\includegraphics[scale=0.09]{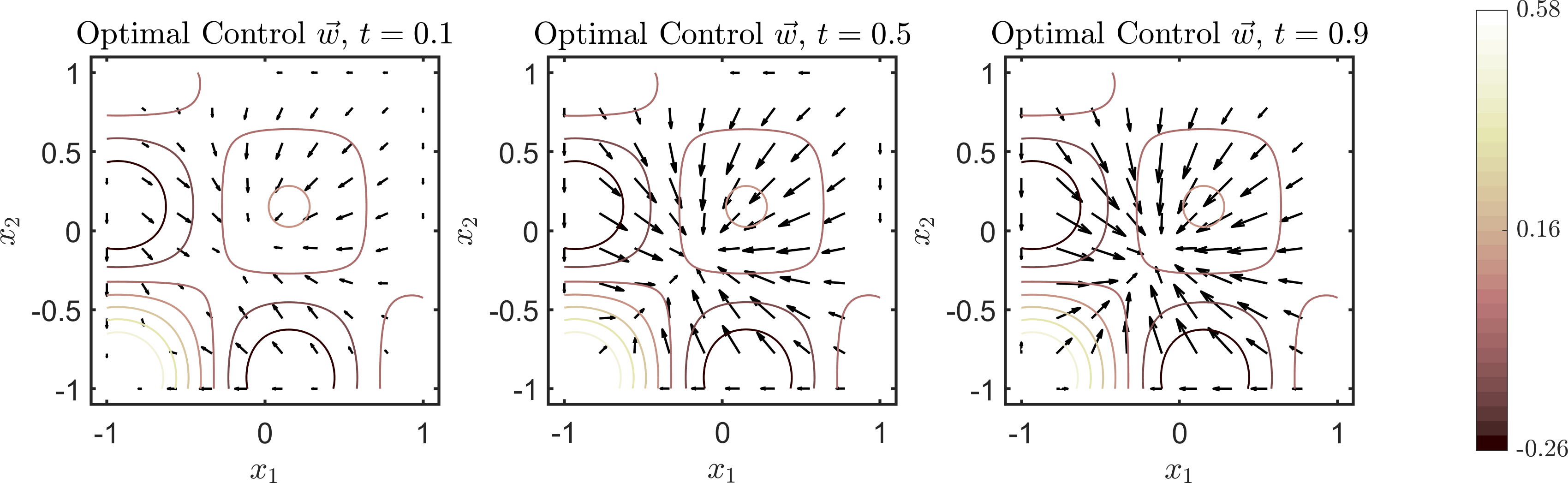}
	\caption{Flow Control, No-Flux: Snapshots of the optimal control for different interaction strengths, $\kappa = -1$, $\kappa = 0$ and $\kappa = 1$ (top to bottom), with $\beta = 10^{-3}$. The lengths of the arrows are proportional to $\|\vec{w}\|$.   A contour plot of the external potential \emph{$V_{\text{ext}}$} is superimposed for reference, with a corresponding colorbar on the right-hand side.} 
	\label{F3ac}
\end{figure}

\input{FCNExample2D.tex}

\begin{figure}[h]
	\centering
	\includegraphics[scale=0.09]{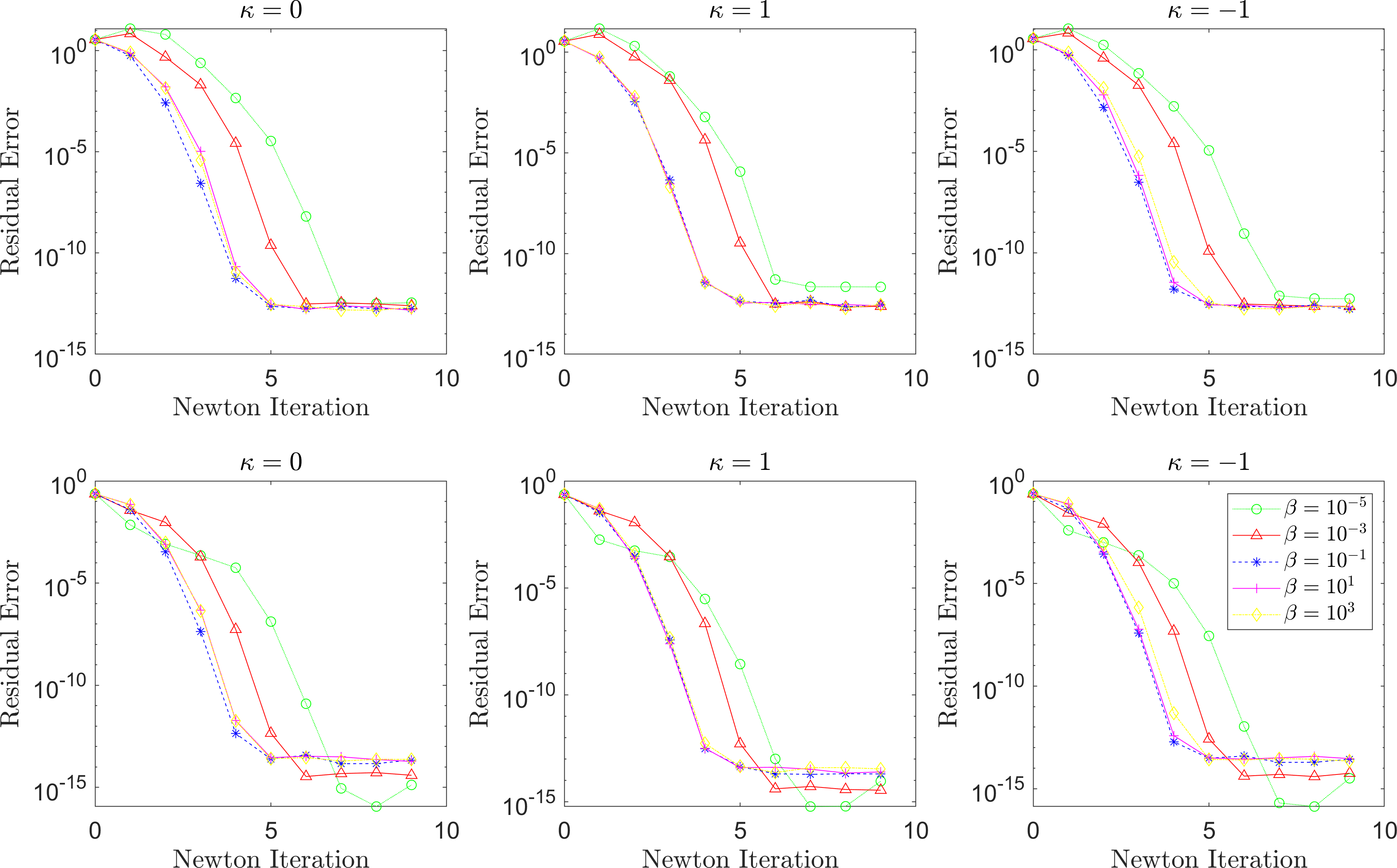}
	\caption{Flow Control, No-Flux: Convergence of the Newton--Krylov algorithm. Top row: convergence in the state variable for different $\kappa$.  Bottom row: convergence in the adjoint variable. Convergence is measured using the residual error.  Note the expected rapid
	convergence in all cases.} 
	\label{Con4}
\end{figure}

 \input{FCNComparison2D}
 
Here we consider \eqref{AdvDiff} with no-flux boundary conditions \eqref{NoFlux}.
The inputs are
\begin{align*}
	\rho_0 &= \frac{1}{4}, \quad \hr = \frac{1}{4}(1-t) + \frac{t}{Z}\exp{\left(-2\left(\left(x_1+0.2\right)^2 + \left(x_2+0.2\right)^2\right)\right)},\\
	V_{\text{ext}} &= \left(\left(x_1 + 0.3\right)^2 - 1\right)\left(\left(x_1-0.4\right)^2 - 0.5\right)
	\left(\left(x_2 + 0.3\right)^2 - 1\right)\left(\left(x_2-0.4\right)^2 - 0.5\right),
\end{align*}
where $Z \approx 1.3791$ is a normalization constant.

In Table  \ref{TabFCN},  the value of the cost functional for the initial configuration ($\mathcal{J}_{uc}$), where $\vec{w} = \vec 0$, is compared with the optimized case ($\mathcal{J}_{c}$) for different values of $\beta$ and for each of the interaction strengths. As expected, in all cases 
$\mathcal{J}_{c} \leq \mathcal{J}_{uc}$  and the lowest values of $\mathcal{J}_{c}$ occur for the smallest $\beta$ values. For large values of $\beta$, applying control is heavily penalized and the optimal control approaches zero, which coincides with the uncontrolled case. 
The numerical solution takes between $200$ and $500$ seconds for the Newton--Krylov solver, and between $500$ and $11500$ seconds for the fixed-point solver.

The results (from the Newton--Krylov scheme) for $\beta = 10^{-3}$ and various interaction strengths, $\kappa$, are shown in Figures \ref{F3a} and  \ref{F3ac}, which display the optimal states and controls, respectively. At earlier times, the density accumulates in regions with potential wells and the areas where the potential is large are avoided. It is clear that the control acts to drive the particle distribution towards the desired state. However, it does not act uniformly around the peak of the desired state, but rather acts strongly in the area between the location of the desired peak and the point $(-1,1)$. This is due to the external potential being large in this area, which requires more control to overcome. It is also evident that more control has to be applied to the repulsive particles, since the desired state requires the particles to accumulate in one part of the domain. In the attractive configuration, the effect of the attraction supports the control action, so less control is needed to reach the desired state. 

Figure \ref{Con4} shows the convergence plot for the residuals arising from the Newton--Krylov scheme for different values of $\beta$. For all values of $\beta$, a residual error of $10^{-13}$ is reached for both state and adjoint variables within $8$ iterations. For larger values of $\beta$, the convergence is slightly faster.

Table \ref{TabFCNCompare} shows the difference $\mathcal{E}$ as defined in Section \ref{sec:Method_Validation} between the solutions of the Newton--Krylov and fixed-point methods. Here we use  $N_1 = N_2 = 30$ and $n = 21$ to ensure accurate solutions with the fixed-point method. As expected, the resulting plots of the density and control for the fixed-point method are very similar to those for the Newton--Krylov solver, and hence we do not show them here. Note that, when the non-local interaction term is turned off ($\kappa = 0$), the error is of order $10^{-3}$ or better, which is reflective of the tolerance chosen for the fixed-point solver. Turning on the interaction term clearly introduces more complexity, and the two solvers disagree slightly more in terms of the optimized state and control. 



\subsubsection{Non-linear (flow) control problem with Dirichlet boundary conditions}

\begin{figure}[h]
	\centering
	\includegraphics[scale=0.09]{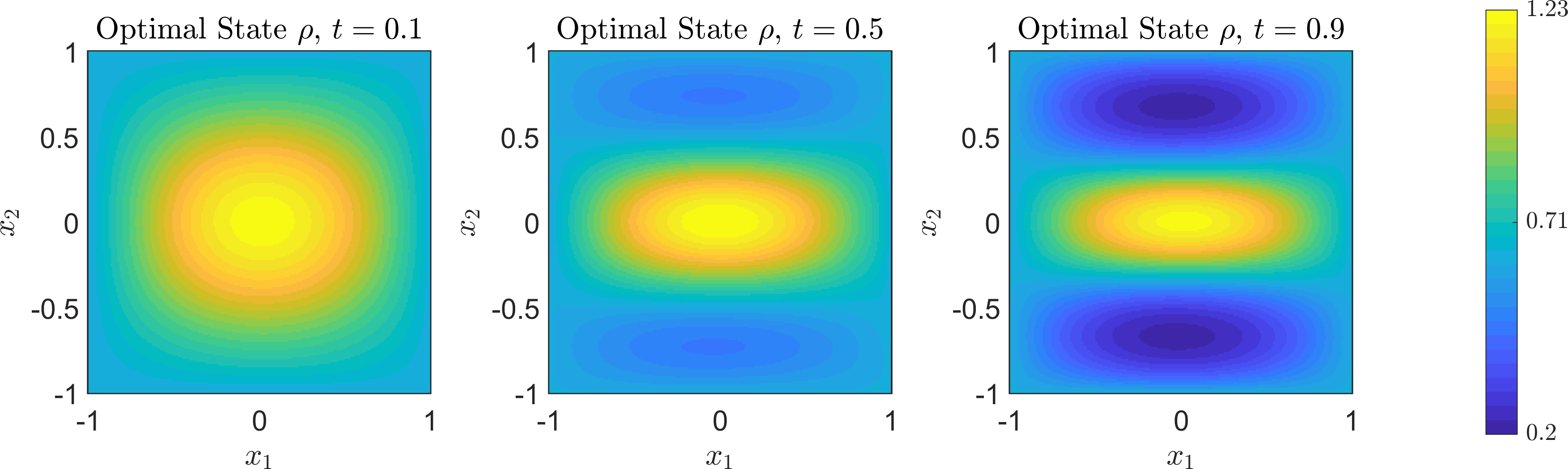}
	\caption{Flow Control, Dirichlet: Snapshots of the optimal $\rho$ for $\kappa = -1$ and $\beta = 10^{-3}$. Note that the optimal states for $\kappa = 0$ and $\kappa = 1$ look almost identical, due to the choice of $\beta$ allowing the control to drive the state close to the desired state $\hr$, and are hence not shown here.} 
	\label{F5a}
\end{figure}

\begin{figure}[h]
	\centering
	\includegraphics[scale=0.09]{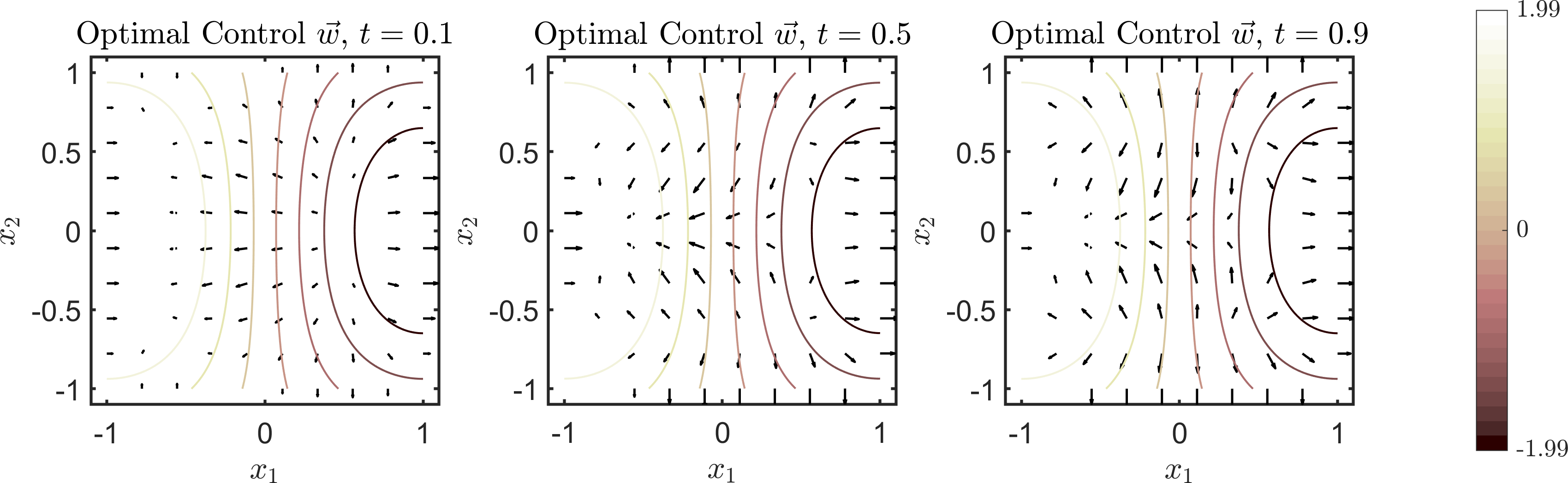}
	\includegraphics[scale=0.09]{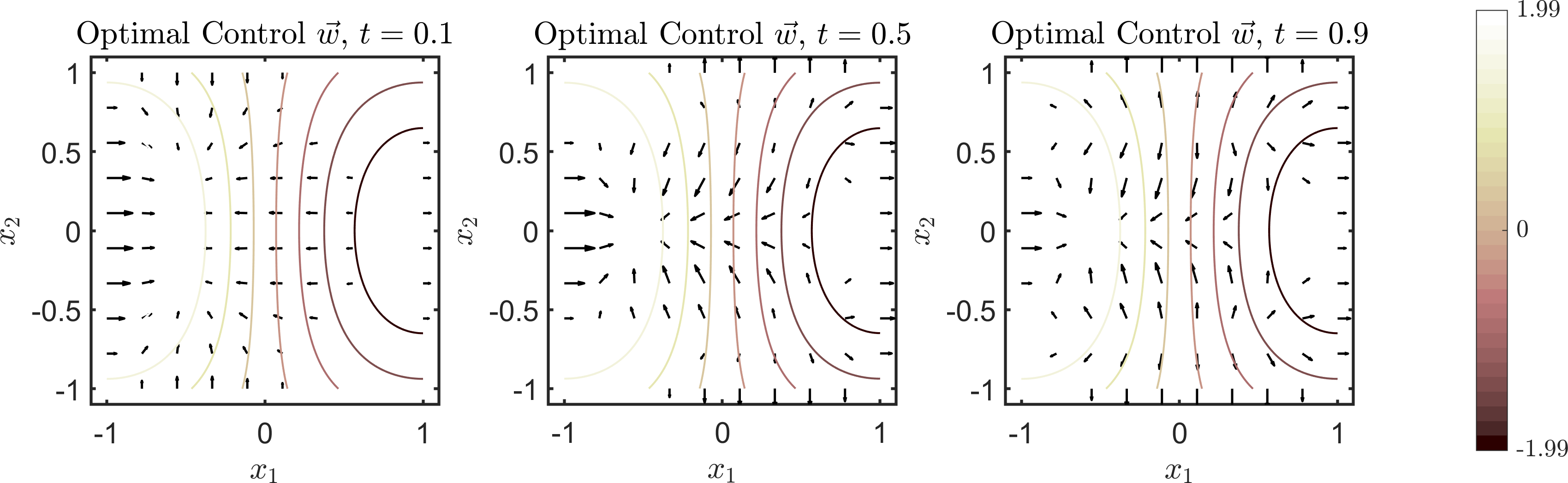}
	\includegraphics[scale=0.09]{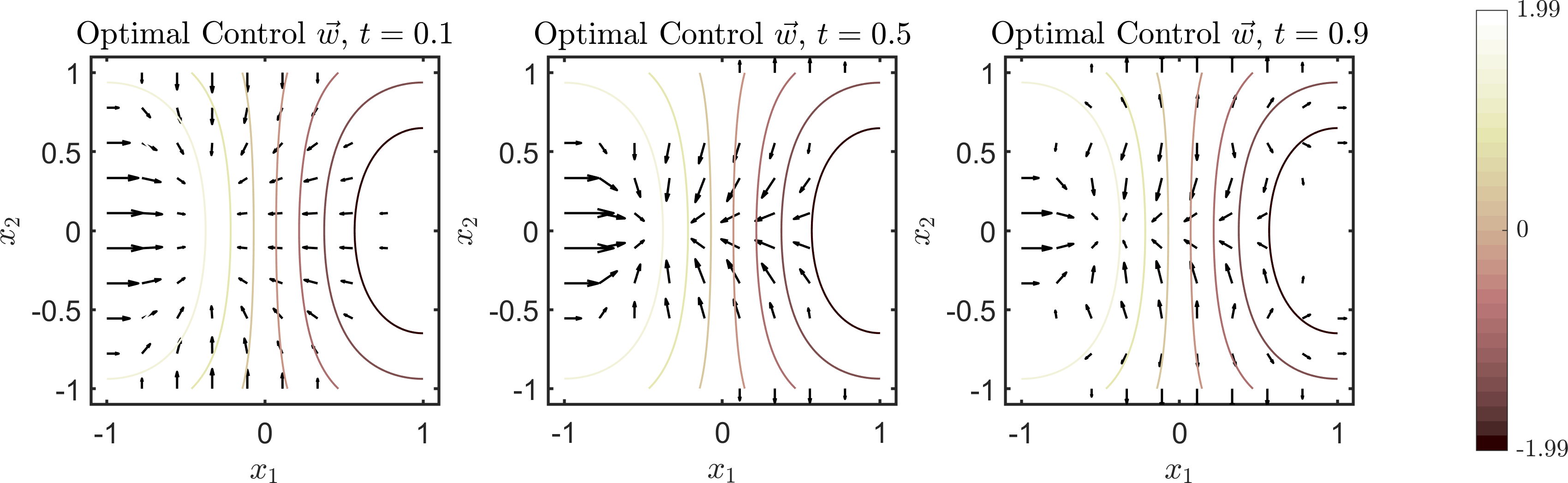}
	\caption{Flow Control, Dirichlet: Snapshots of the optimal control for $\kappa = -1$, $\kappa = 0$ and $\kappa = 1$ (top to bottom), 
	for $\beta = 10^{-3}$. See Figure~\ref{F3ac} for further details.}
	\label{F5ac}
\end{figure}

\input{FCDExample2D.tex}

Our next example is a control problem of type \eqref{AdvDiff}, with Dirichlet boundary conditions \eqref{Dirichlet}, and
\begin{align*}
	\rho_0 &= \left(\frac{\pi}{4}\right)^2\cos\left(\frac{\pi x_1}{2}\right)\cos\left(\frac{\pi x_2}{2}\right) + \left(\frac{\pi}{4}\right)^2, \quad 
	V_{\text{ext}} = 2\sin\left(\frac{\pi x_1}{2}\right) \sin\left(\frac{\pi x_2}{3} - \frac{\pi}{2}\right),\\
	\hr &= (1 - t)\left(\left(\frac{\pi}{4}\right)^2\cos\left(\frac{\pi x_1}{2}\right)\cos\left(\frac{\pi x_2}{2}\right) + \left(\frac{\pi}{4}\right)^2\right) + t\left(\left(\frac{\pi}{4}\right)^2\cos\left(\frac{\pi x_1}{2}\right)\cos\left(\frac{3\pi x_2}{2}\right) + \left(\frac{\pi}{4}\right)^2\right).
\end{align*}
The resulting costs are in Table \ref{TabFCD}. The solution of each example takes between $100$ and $200$ seconds.

The desired state prescribes the density to move from one uniform bump in the middle of the domain, to accumulate in a steeper, elongated shape across the $x_1$-axis. In this example, the effect of the different interaction strengths and the external potential on the control is clearly visible; see Figure \ref{F5ac}. The external potential is large on the left side of the domain, which naturally drives the density away from this region, so that most effort of the control variable is concentrated on the right side of the domain. 
However, for attractive particles, additional control must be applied in the right side of the domain, since the attractive particles oppose the density spreading out along the $x_1$-axis. For these attractive particles, the initially clumped density is a more `natural' state; see Figure \ref{F5a}. In contrast, for repulsive particles, most of the work of the control is done to push the particles together. These controls can be seen in Figure \ref{F5ac}.


\subsubsection{Linear (source) control problem with no-flux boundary conditions}

\begin{figure}[h]
	\centering
	\includegraphics[scale=0.09]{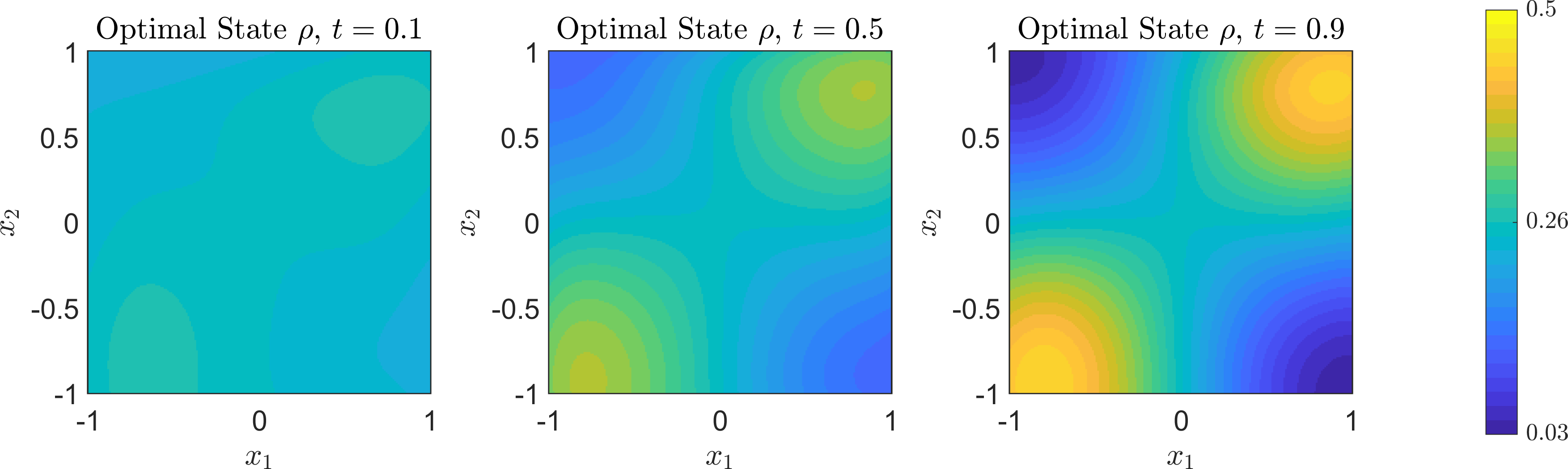}
	\includegraphics[scale=0.09]{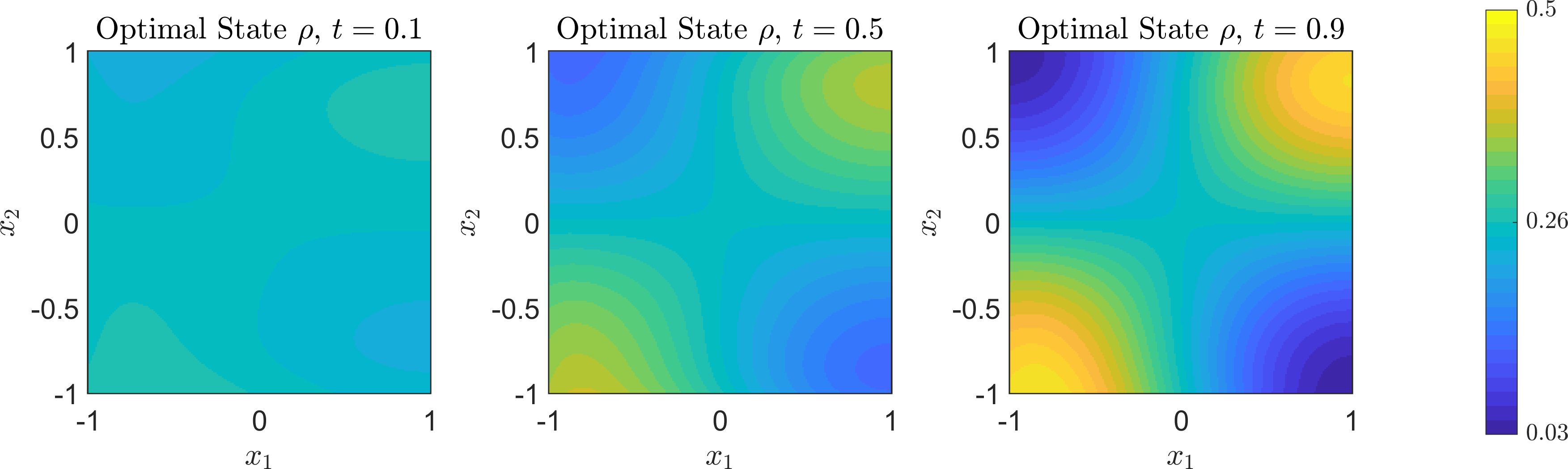}	
	\includegraphics[scale=0.09]{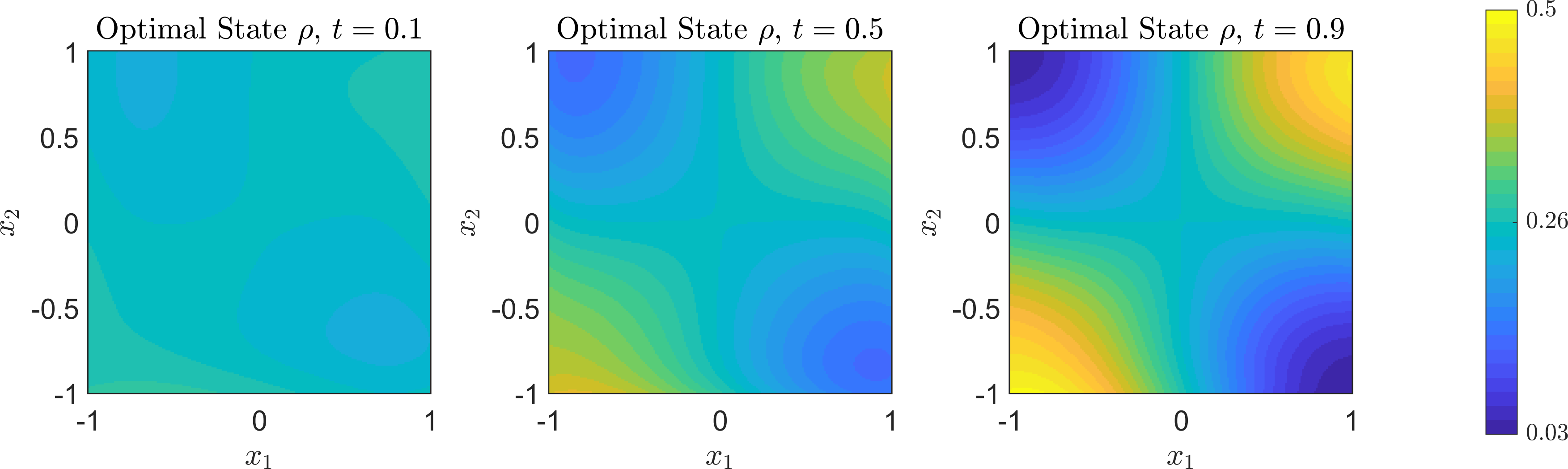}
	\caption{Source Control, No-Flux: Snapshots of the optimal $\rho$ for $\kappa = -1$, $\kappa = 0$, and $\kappa = 1$
	(top to bottom), for $\beta = 10^{-3}$.} 
	\label{FSCN1}
\end{figure}


\input{SCNExample2D.tex}

We consider problem \eqref{AdvDiff_Linear} with no-flux boundary conditions \eqref{NoFlux_Linear}.
The chosen inputs for this example are
\begin{align*}
	\rho_0 &= \frac{1}{4}, \quad V_{\text{ext}} = \cos\left(\frac{\pi x_1}{5} - \frac{\pi}{5}\right)\sin\left(\frac{\pi x_2}{5}\right),\\
	\hr &= \frac{1}{4}(1 - t) + t\left(\frac{1}{4}\sin\left(\frac{\pi \left(x_1 - 2\right)}{2}\right)\sin\left(\frac{\pi (x_2 - 2)}{2}\right) + \frac{1}{4}\right).
\end{align*}
The resulting costs for different $\beta$ and $\kappa$ are in Table  \ref{TabSCN}. We note that the solution of each example takes between $100$ and $200$ seconds, apart from the $\beta = 10^{-5}$ case which takes around $25$ seconds.

In Figure \ref{FSCN1} we show (for $\beta = 10^{-3}$) the optimal states for different interaction strengths. 
Since $\beta$ is small, the optimal state is very close to the desired state $\hr$ (not shown). We can observe clear effects on the optimal state and the control from the external potential $V_{\text{ext}}$. Since $V_{\text{ext}}$ is large around $x_2 = 1$, more control has to be applied in this area to force the density towards $\widehat{\rho}$. 
It can also be seen that the state is slightly asymmetric because of this effect, despite $\hr$ being symmetric.

The effect of the different interaction strengths on the state can be observed in Figure \ref{FSCN1}, by inspecting the shape of the particle distribution. The desired state $\hr$ prescribes higher density near the two corners $(-1,1)$ and $(1,1)$. Without control or an external potential, repulsive particles accumulate on the boundary of the domain, whilst attractive particles favour the centre of the domain.  Hence in this example, where the target density is higher near the boundary, less control needs to be applied for repulsive particles. 
For attractive particles, the accumulated particles are arranged in a rounder shape, while the repulsive particles are more spread out, as would be expected from their interactions.


\subsubsection{Linear (source) control problem with Dirichlet boundary conditions}

\begin{figure}[h]
	\centering
	\includegraphics[scale=0.09]{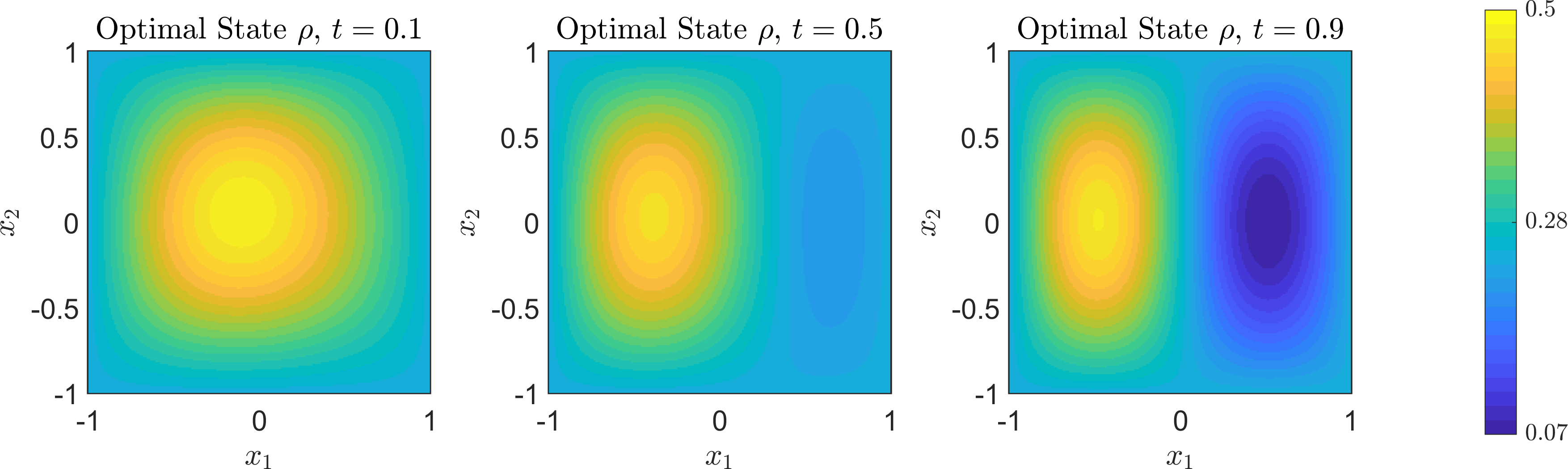}
	\caption{Source Control, Dirichlet: Snapshots of the optimal $\rho$ for $\kappa = -1$, for $\beta = 10^{-3}$. Note that the optimal states for $\kappa = 0$ and $\kappa = 1$ look identical, due to the choice of $\beta$ allowing the control to drive the state close to $\hr$.}
	\label{F2a}
\end{figure}

\begin{figure}[h]
	\centering
	\includegraphics[scale=0.09]{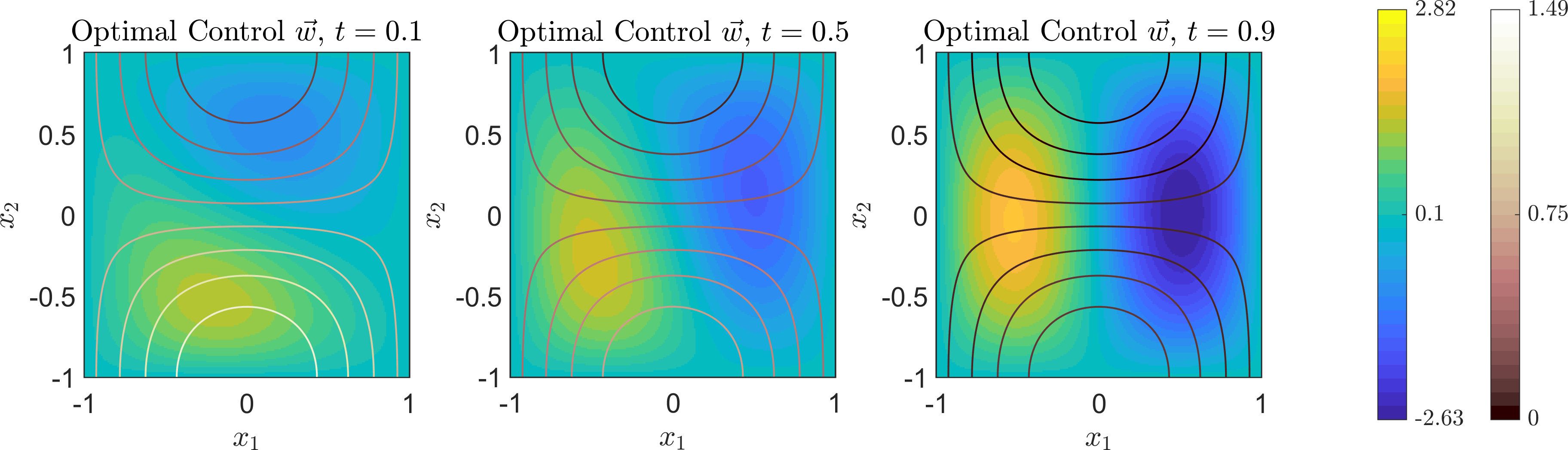}
	\includegraphics[scale=0.09]{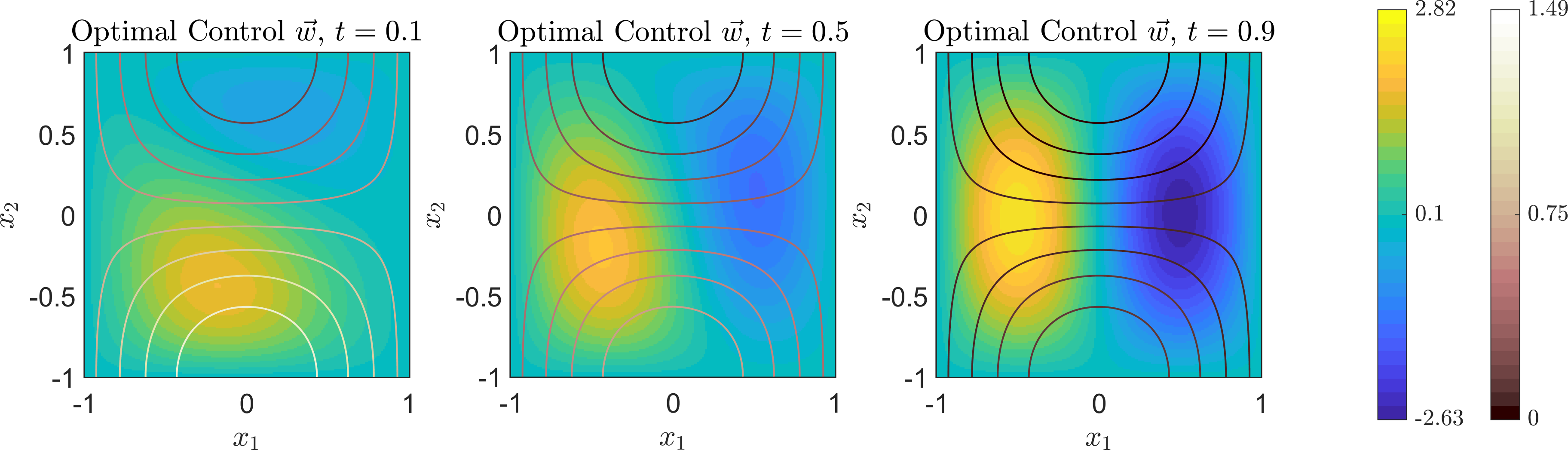}
	\includegraphics[scale=0.09]{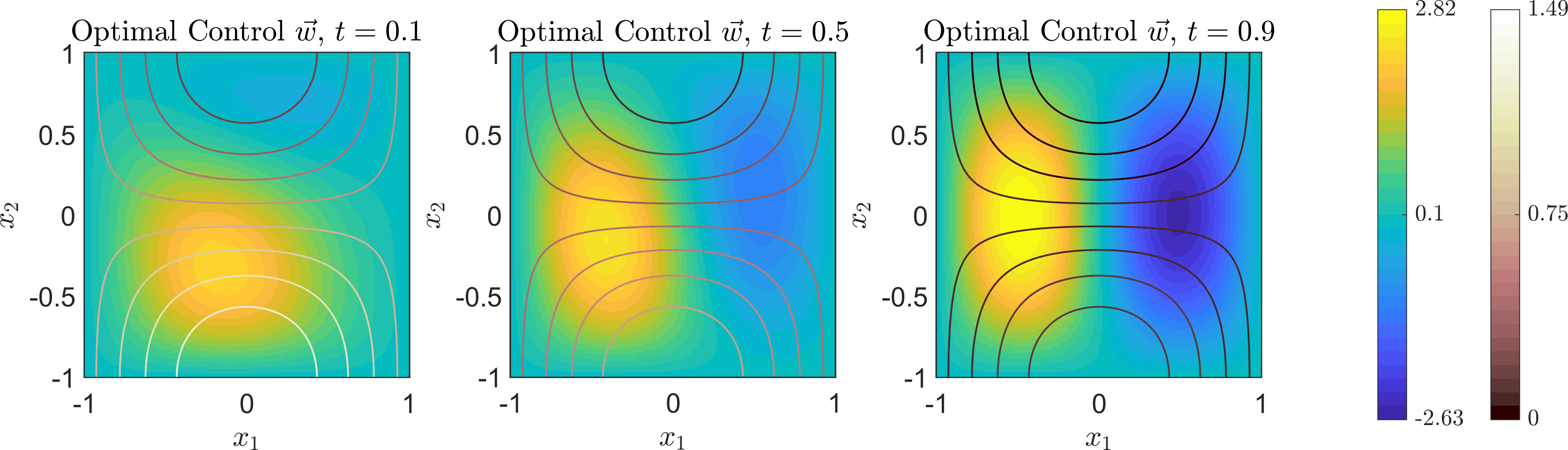}
	\caption{Source Control, Dirichlet: Snapshots of the optimal control for $\kappa = -1$, $\kappa = 0$, and $\kappa = 1$ 
	(top to bottom), with $\beta = 10^{-3}$. A contour plot of the external potential \emph{$V_{\text{ext}}$} is superimposed on the control plots for reference, with a corresponding colorbar on the right-hand side.}
	\label{F2ac}
\end{figure}

\input{SCDExample2D.tex}

We consider problem \eqref{AdvDiff_Linear} with Dirichlet boundary conditions \eqref{Dirichlet}.
The chosen inputs are
\begin{align*}
	\rho_0 &= \frac{1}{4}\cos\left(\frac{\pi x_1}{2}\right)\cos\left(\frac{\pi x_2}{2}\right) + \frac{1}{4}, \quad V_{\text{ext}} =  \frac{3}{4}(1-t)\left(-\cos\left(\frac{\pi x_1}{2}\right)\sin\left(\frac{\pi x_2}{2}\right) + 1\right),\\
	\hr &= (1 - t)\left(\frac{1}{4}\cos\left(\frac{\pi x_1}{2}\right)\cos\left(\frac{\pi x_2}{2}\right) + \frac{1}{4}\right) - t\left(\frac{1}{4}\sin\left(\pi x_1\right)\sin\left(\frac{\pi x_2}{2} - \frac{\pi}{2}\right) + \frac{1}{4}\right).
\end{align*}
Note, in particular, that the  external potential is time dependent. Since it decays over time, this results in the strongest effect of $V_{\text{ext}}$ being visible at earlier times.  The resulting costs for different $\beta$ and $\kappa$ can be seen in Table \ref{TabSCD}.  The solution of each example takes between $70$ and $100$ seconds, apart from the $\beta = 10^{-5}$ case which takes around $25$ seconds.

We show the optimal state for $\beta = 10^{-3}$ and varying $\kappa$ in Figure \ref{F2a}, which is once again close to the target state, $\hr$, irrespective of the interaction strength.  The corresponding optimal controls are shown in Figure \ref{F2ac}.
Since the external potential is large around the bottom half of the domain, the density is not centred in the middle of the domain, but shifted slightly upwards. At the same time it can be observed that at $t = 0.1$, the control is applied where the external potential is steep. At later times, the control is mostly applied where the density is prescribed to accumulate approximately in the form of the desired state $\hr$, which is at the left half of the domain.
While the qualitative behaviour of the control is similar in each case, it can be seen that less control has to be applied for attractive particles compared to repulsive ones, since the attraction causes the particles to clump together, which supports the shape of the desired state $\hr$.


\subsection{Three-Dimensional Example}


	\begin{figure}[h]
		\centering
		\includegraphics[scale=0.09]{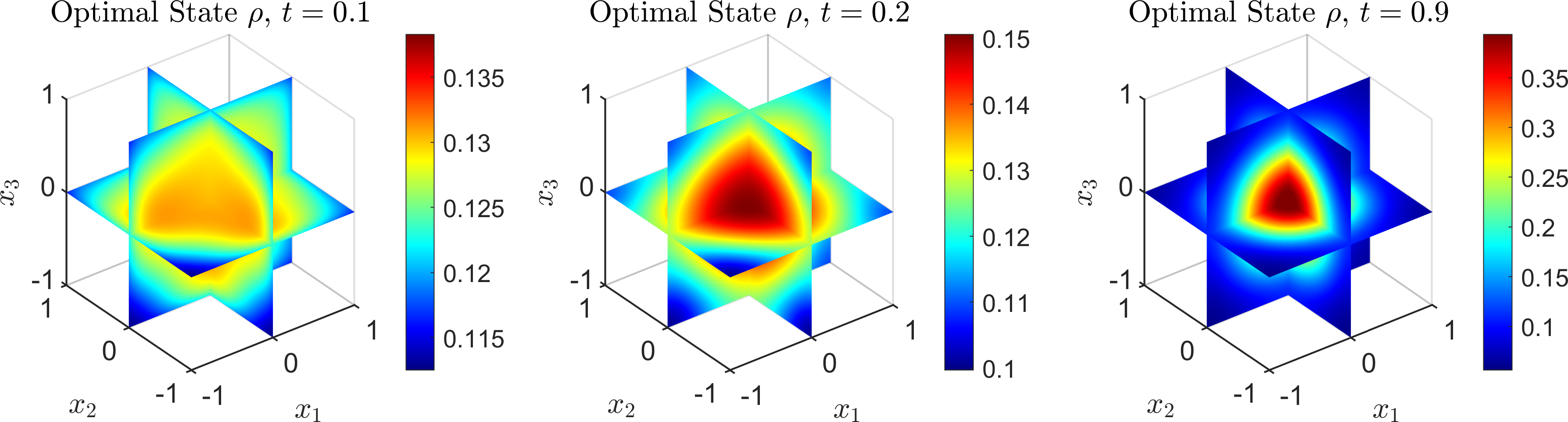}
		\includegraphics[scale=0.09]{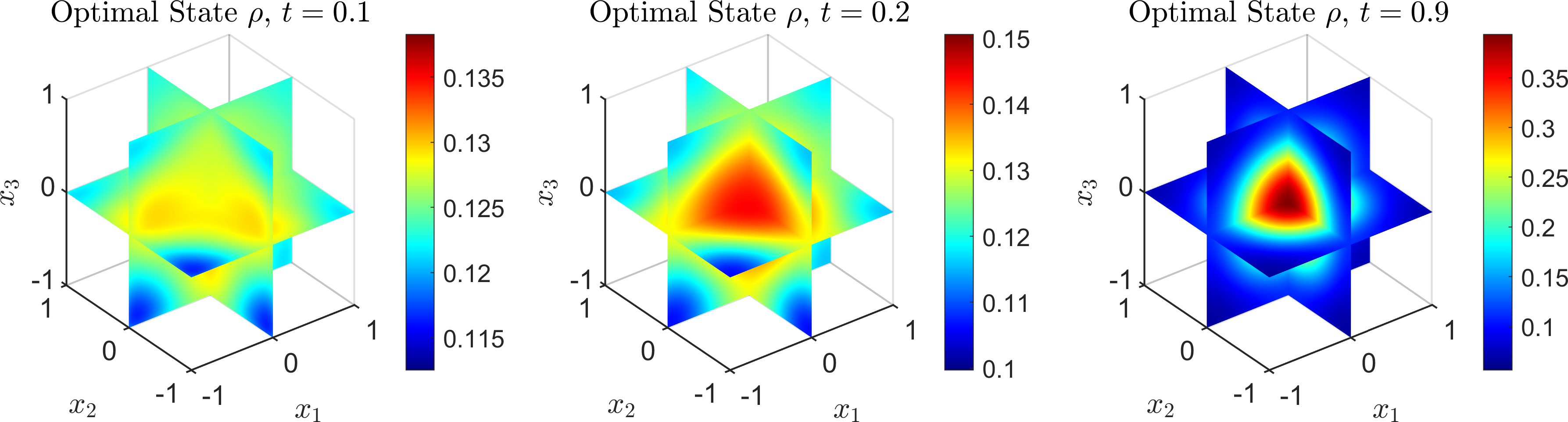}
		\includegraphics[scale=0.09]{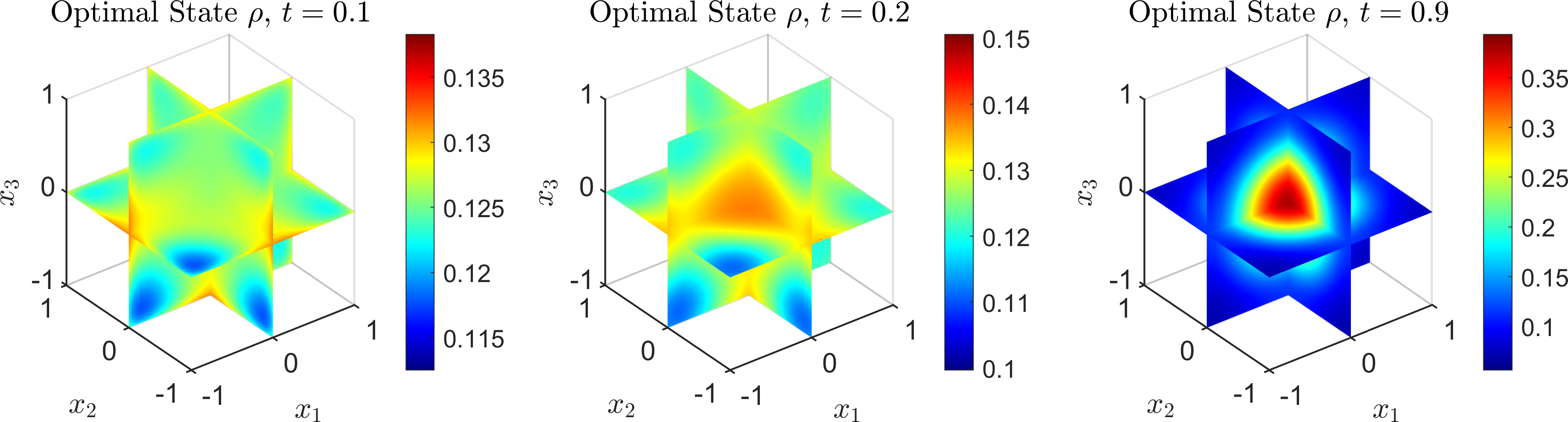}
		\caption{3D Flow Control, No-Flux: Snapshots of the optimal state $\rho$ for $\kappa = -1$, $\kappa = 0$ and $\kappa = 1$
		(top to bottom), for $\beta = 10^{-3}$.} 
		\label{F2}
	\end{figure}


Our final example is three-dimensional, featuring the non-linear flow control problem \eqref{AdvDiff} with no-flux boundary conditions \eqref{NoFlux}.  This has been chosen as an illustrative example in three dimensions since it is both the most challenging combination of control type and boundary conditions and also the most physically relevant for our applications.
The chosen inputs are
\begin{align*}
	\rho_0 ={}& \frac{1}{8}, \quad \hr = \frac{1}{8}(1-t) + t\left(\frac{\pi}{4}\right)^3\cos\left(\frac{\pi x_1}{2}\right)\cos\left(\frac{\pi x_2}{2}\right)\cos\left(\frac{\pi x_3}{2}\right),\\
	V_{\text{ext}} ={}& \left(\left(x_1 + 0.3\right)^2 - 1\right)\left(\left(x_1 - 0.4\right)^2 - 0.5\right) \left(\left(x_2 + 0.3\right)^2 - 1\right)\left(\left(x_2 - 0.4\right)^2 - 0.5\right)\\
	& \quad\left(\left(x_3 + 0.3\right)^2 - 1\right)\left(\left(x_3 - 0.4\right)^2 - 0.5\right).
\end{align*}
This example is only run for $\beta = 10^{-3}$, due to a running time of approximately $35$ hours per problem.  This is a simple consequence
of the `curse of dimensionality'.

The effect of the different interaction strengths is clearly displayed in Figure \ref{F2} and is particularly obvious in earlier times of the particle evolution. 
It is evident that attractive particles enhance the control in pushing the density into a cluster in the middle of the domain, as prescribed by the desired state $\hr$, while more control is needed for a similar effect in the repulsive setup.
We get, for $\kappa = 0$, $\mathcal J_c = 0.0078$. This can be compared to $\mathcal J_{uc} = 0.0195$ from the computed forward problem with $\vec{w} = \vec 0$. For $\kappa = 1$ we obtain $\mathcal J_c = 0.0102$, compared to $\mathcal J_{uc} = 0.0232$ in the uncontrolled case, and for $\kappa = -1$ we have $\mathcal J_c = 0.0059$, with $\mathcal J_{uc} = 0.0477$.  As expected, the optimal control leads to a cost
which is significantly lower than in the uncontrolled case.

\section{\label{sec:Conc} Concluding Remarks}
We have derived an accurate and efficient algorithmic strategy for solving the first-order optimality conditions arising from PDE-constrained optimization problems, along with additional integral terms, describing multiscale particle dynamics problems. Our approach, linked to the DDFT approach applied to (non-optimized) systems in statistical mechanics, applies a pseudospectral method in space and time, and utilizes fixed-point and Newton--Krylov schemes within the optimization solver. This novel methodology is more general in scope than existing numerical implementations for similar problems, and exhibits the substantial computational benefits of applying such methods for non-local, non-linear systems of PDEs. Numerical tests indicate the potency of our approach for a range of examples, boundary conditions, and problem parameters. An open-source software implementation of our methodology is available at~\cite{2DChebClassPDECO}. There are many possible extensions to our approach: for instance, one may apply our methodology to problems where the misfit between state and desired state is measured at some final time only, models with different cost functionals, and boundary control problems. Furthermore, methods of this type may be tailored to specific particle dynamics applications, in fields such as opinion dynamics, flocking, swarming, and optimal control problems in robotics, and such applications will be tackled in future work.

\section*{Acknowledgements}
MA and JCR are supported by The Maxwell Institute Graduate School in Analysis and its Applications 
(EPSRC grant EP/L016508/01), the Scottish Funding Council, Heriot-Watt University, and The University of Edinburgh.
BDG acknowledges support from EPSRC grant EP/L025159/1;
JWP from EPSRC grant EP/S027785/1 and an Alan Turing Institute Fellowship. 

\appendix

\section{Validation of Newton--Krylov and Fixed-Point Methods} \label{appendix:validation}
We compare results from the Newton--Krylov and fixed-point algorithms to establish that both methods correctly solve a problem with a known exact solution.  The methods are also compared  for the flow control problem with no-flux boundary conditions in
Section~\ref{sec:fcn}; see Table~\ref{TabFCNCompare} in particular.  The example here is also 
of this form; the differences here are that there are no interactions, and we include an additional source term such that the problem has an analytic solution.

The exact (analytic) solutions and input choices are
\begin{align*}
\rho_{ex} ={}& \frac{1}{4}\beta^{1/2} e^t \left(\cos(\pi x_1) +1\right)\left(\cos(\pi x_2) + 1\right), \\
\adj_{ex} ={}& \frac{1}{4}\beta^{1/2}\left(e^T - e^t\right) \left(\cos(\pi x_1) +1\right)\left(\cos(\pi x_2) + 1\right),\\
\vec{w}_{ex} ={}& \frac{\pi}{16} e^t\left(e^T - e^t\right) \left(\cos(\pi x_1) +1\right)\left(\cos(\pi x_2) + 1\right)\left[\begin{array}{c}
\sin(\pi x_1) \left(\cos(\pi x_2) + 1\right)\\
\left(\cos(\pi x_1) +1\right)\sin(\pi x_2) 
\end{array}\right],\\
V_{\text{ext}} ={}& \cos(\pi x_1) \cos(\pi x_2),\\
\hr ={}& -  \frac{\pi^2}{4}\beta^{1/2}\left(e^T - e^t\right)\left[
\cos(\pi x_1) \left(\cos(\pi x_2) + 1\right)
 + \cos(\pi x_2) \left(\cos(\pi x_1) + 1\right) \right.\\
 & \quad \left. + \sin^2\left(\pi x_1\right)\cos(\pi x_2)\left(\cos(\pi x_2)+1\right)  
 + \sin^2\left(\pi x_2\right)\cos(\pi x_1)\left(\cos(\pi x_1)+1\right)  \right]\\
& \quad + \frac{\pi^2}{2}\beta^{1/2} e^t \left(e^T - e^t\right)^2
\cos^4 \left(\frac{\pi x_1}{2}\right)\cos^4 \left(\frac{\pi x_2}{2}\right)\left(\cos(\pi x_1) \cos(\pi x_2) - 1\right),\\
f ={}& \frac{1}{4}\beta^{1/2} e^t \left(\cos(\pi x_1) +1\right)\left(\cos(\pi x_2) + 1\right) + \frac{\pi^2}{4}\beta^{1/2} e^t \left[\cos( \pi x_1) \left(2 \cos(\pi x_2) + 1 \right) + \cos(\pi x_2) \right]\\
&\quad - \pi^2\beta^{1/2}e^{t} \cos^2 \left(\frac{\pi x_1}{2}\right)\cos^2 \left(\frac{\pi x_2}{2}\right)\left[\cos(\pi x_1)\left(1 - 4 \cos(\pi x_2)\right) + \cos(\pi x_2)\right]\\
& \quad + \frac{\pi^2}{4}\beta^{1/2}e^{2t}\left(e^T - e^t\right) \cos^4\left(\frac{\pi x_1}{2}\right) \cos^4\left(\frac{\pi x_2}{2}\right)\left[\cos(\pi x_1) \left(6 \cos(\pi x_2) + 1\right) + \cos(\pi x_2) - 4 \right].
\end{align*}

The results for this example can be seen in Table \ref{TabFCNExact}, which displays the error of the computed state and adjoint variables using the Newton--Krylov solver in reference to the exact solution after $10$ iterations, with $n = 10$ and $N_1 = N_2 = 20$ points, and the fixed-point solver after $1$ iteration having been supplied with the exact solution for $\vec{w}$ as an initial guess, when using $n = 21$ and $N_1 = N_2 = 30$ points. Each example takes between $38$ and $88$ seconds to run for the fixed-point solver, while it takes $80$ seconds to solve the exact problem with the Newton--Krylov algorithm. It should be noted however that the Newton--Krylov solver requires an initial guess for $\rho$ and $\adj$ at all time points, while the fixed-point solver requires an initial guess for $\vec w$ at all times. The Newton--Krylov initial guess at all time points consists of the initial and final time conditions $\rho_0$ and $\adj_T$, 
while the initial guess for the fixed-point algorithm is $\vec w_{ex}$. 
The errors $\mathcal{E}$ (as defined in Section \ref{sec:Method_Validation}) displayed in Table \ref{TabFCNExact}  are small for both methods.
However, the errors made by the Newton--Krylov solver are several orders smaller than those of the fixed-point algorithm, demonstrating the superiority of the higher-order method over the first-order method. It is therefore the natural choice for solving the optimization problems considered in this paper, although the accuracy 
achieved by the fixed-point scheme is sufficient to retrieve good results in cases where the Newton--Krylov scheme cannot easily be applied, such as in optimal control problems with box constraints.

\input{FCNExact2D}

\bibliographystyle{plain}
\bibliography{ParticleDynamicsPDECO}

\end{document}

%% file: FCNExample2D.tex
\begin{table}
\centering
\begin{small}
\begin{tabular}{ | c | c || c | c | c | c | c |}
\hline
\multicolumn{2}{|c||}{}& $\beta = 10^{-5}$ & $\beta = 10^{-3}$ & $\beta = 10^{-1}$ & $\beta = 10^{1}$ & $\beta = 10^{3}$  \\
\hline
\hline
\multirow{2}{*}{$\kappa= \numprint{0}$}  & $\mathcal{J}_{uc}$ & $\numprint{2.67e-2}$ & $\numprint{2.67e-2}$ & $\numprint{2.67e-2}$ & $\numprint{2.67e-2}$ & $\numprint{2.67e-2}$\\
 & $\mathcal{J}_c$ & $\numprint{8.23e-5}$ & $\numprint{3.87e-3}$ & $\numprint{2.50e-2}$ & $\numprint{2.67e-2}$ & $\numprint{2.67e-2}$\\
\hline
\multirow{2}{*}{$\kappa= \numprint{1}$}  & $\mathcal{J}_{uc}$ & $\numprint{3.29e-2}$ & $\numprint{3.29e-2}$ & $\numprint{3.29e-2}$ & $\numprint{3.29e-2}$ & $\numprint{3.29e-2}$\\
 & $\mathcal{J}_c$ & $\numprint{1.16e-4}$ & $\numprint{5.44e-3}$ & $\numprint{3.13e-2}$ & $\numprint{3.29e-2}$ & $\numprint{3.29e-2}$\\
\hline
\multirow{2}{*}{$\kappa= \numprint{-1}$}  & $\mathcal{J}_{uc}$ & $\numprint{2.09e-2}$ & $\numprint{2.09e-2}$ & $\numprint{2.09e-2}$ & $\numprint{2.09e-2}$ & $\numprint{2.09e-2}$\\
 & $\mathcal{J}_c$ & $\numprint{5.71e-5}$ & $\numprint{2.63e-3}$ & $\numprint{1.92e-2}$ & $\numprint{2.09e-2}$ & $\numprint{2.09e-2}$\\
\hline
\end{tabular}
\end{small}
\caption{Flow Control, No-Flux: Cost when $\vec{w}=\vec{0}$ ($\mathcal{J}_{uc}$) and optimal control cost ($\mathcal{J}_{c}$) for a range of $\kappa$ and $\beta$.}
\label{TabFCN}
\end{table}

%% file: FCNComparison2D.tex
\begin{table}
\centering
\begin{small}
\begin{tabular}{ | c | c || c | c | c | c | c |}
\hline
\multicolumn{2}{|c||}{}& $\beta = 10^{-3}$ & $\beta = 10^{-2}$ & $\beta = 10^{-1}$ & $\beta = 10^{1}$ & $\beta = 10^{3}$  \\
\hline
\hline
\multirow{2}{*}{$\kappa= \numprint{0}$}  & $\mathcal{E}_{\rho}$ & $\numprint{2.48e-4}$ & $\numprint{2.34e-4}$ & $\numprint{2.34e-4}$ & $\numprint{2.34e-4}$ & $\numprint{2.34e-4}$\\
 & $\mathcal{E}_{\vec{w}}$ & $\numprint{9.69e-3}$ & $\numprint{9.64e-4}$ & $\numprint{9.61e-5}$ & $\numprint{9.61e-7}$ & $\numprint{9.61e-9}$\\
\hline
\multirow{2}{*}{$\kappa= \numprint{1}$}  & $\mathcal{E}_{\rho}$ & $\numprint{3.66e-2}$ & $\numprint{2.48e-2}$ & $\numprint{1.85e-2}$ & $\numprint{2.00e-2}$ & $\numprint{2.01e-2}$\\
 & $\mathcal{E}_{\vec{w}}$ & $\numprint{5.53e-2}$ & $\numprint{3.47e-2}$ & $\numprint{9.63e-3}$ & $\numprint{1.17e-4}$ & $\numprint{1.17e-6}$\\
\hline
\multirow{2}{*}{$\kappa= \numprint{-1}$}  & $\mathcal{E}_{\rho}$ & $\numprint{5.38e-2}$ & $\numprint{8.90e-2}$ & $\numprint{5.73e-2}$ & $\numprint{4.73e-2}$ & $\numprint{4.72e-2}$\\
 & $\mathcal{E}_{\vec{w}}$ & $\numprint{7.95e-2}$ & $\numprint{6.19e-2}$ & $\numprint{1.11e-2}$ & $\numprint{1.44e-4}$ & $\numprint{1.44e-6}$\\
\hline
\end{tabular}
\end{small}
\caption{Flow Control, No-Flux: Comparison of $\rho$ and $\vec{w}$ obtained from Newton--Krylov and fixed-point--Armijo--Wolfe solvers, for a range of $\kappa$ and $\beta$. }
\label{TabFCNCompare}
\end{table}

%% file: FCDExample2D.tex
\begin{table}
\centering
\begin{small}
\begin{tabular}{ | c | c || c | c | c | c | c |}
\hline
\multicolumn{2}{|c||}{}& $\beta = 10^{-5}$ & $\beta = 10^{-3}$ & $\beta = 10^{-1}$ & $\beta = 10^{1}$ & $\beta = 10^{3}$  \\
\hline
\hline
\multirow{2}{*}{$\kappa= \numprint{0}$}  & $\mathcal{J}_{uc}$ & $\numprint{1.58e-1}$ & $\numprint{1.58e-1}$ & $\numprint{1.58e-1}$ & $\numprint{1.58e-1}$ & $\numprint{1.58e-1}$\\
 & $\mathcal{J}_c$ & $\numprint{4.15e-4}$ & $\numprint{7.74e-3}$ & $\numprint{1.30e-1}$ & $\numprint{1.58e-1}$ & $\numprint{1.58e-1}$\\
\hline
\multirow{2}{*}{$\kappa= \numprint{1}$}  & $\mathcal{J}_{uc}$ & $\numprint{2.12e-1}$ & $\numprint{2.12e-1}$ & $\numprint{2.12e-1}$ & $\numprint{2.12e-1}$ & $\numprint{2.12e-1}$\\
 & $\mathcal{J}_c$ & $\numprint{4.16e-4}$ & $\numprint{1.03e-2}$ & $\numprint{1.85e-1}$ & $\numprint{2.12e-1}$ & $\numprint{2.12e-1}$\\
\hline
\multirow{2}{*}{$\kappa= \numprint{-1}$}  & $\mathcal{J}_{uc}$ & $\numprint{4.03e-1}$ & $\numprint{4.03e-1}$ & $\numprint{4.03e-1}$ & $\numprint{4.03e-1}$ & $\numprint{4.03e-1}$\\
 & $\mathcal{J}_c$ & $\numprint{4.53e-4}$ & $\numprint{8.65e-3}$ & $\numprint{1.74e-1}$ & $\numprint{3.87e-1}$ & $\numprint{4.03e-1}$\\
\hline
\end{tabular}
\end{small}
\caption{Flow Control, Dirichlet: Cost when $\vec{w}=\vec{0}$ ($\mathcal{J}_{uc}$) and optimal control cost ($\mathcal{J}_{c}$) for a range of $\kappa$ and $\beta$. For $\beta = 10$, the cost functionals differ by $10^{-4}$ for $\kappa = 0$ and $\kappa = 1$, and by $10^{-2}$ for $\kappa = -1$. For $\beta = 10^3$, the cost functionals differ by $10^{-7}$ for $\kappa = 0$, by $10^{-6}$ for $\kappa = 1$, and by $10^{-4}$ for $\kappa = -1$.}
\label{TabFCD}
\end{table}

%% file: SCNExample2D.tex
\begin{table}
\centering
\begin{small}
\begin{tabular}{ | c | c || c | c | c | c | c |}
\hline
\multicolumn{2}{|c||}{}& $\beta = 10^{-5}$ & $\beta = 10^{-3}$ & $\beta = 10^{-1}$ & $\beta = 10^{1}$ & $\beta = 10^{3}$  \\
\hline
\hline
\multirow{2}{*}{$\kappa= \numprint{0}$}  & $\mathcal{J}_{uc}$ & $\numprint{1.90e-2}$ & $\numprint{1.90e-2}$ & $\numprint{1.90e-2}$ & $\numprint{1.90e-2}$ & $\numprint{1.90e-2}$\\
 & $\mathcal{J}_c$ & $\numprint{1.29e-5}$ & $\numprint{6.65e-4}$ & $\numprint{1.37e-2}$ & $\numprint{1.89e-2}$ & $\numprint{1.90e-2}$\\
\hline
\multirow{2}{*}{$\kappa= \numprint{1}$}  & $\mathcal{J}_{uc}$ & $\numprint{1.94e-2}$ & $\numprint{1.94e-2}$ & $\numprint{1.94e-2}$ & $\numprint{1.94e-2}$ & $\numprint{1.94e-2}$\\
 & $\mathcal{J}_c$ & $\numprint{1.59e-5}$ & $\numprint{7.43e-4}$ & $\numprint{1.42e-2}$ & $\numprint{1.93e-2}$ & $\numprint{1.94e-2}$\\
\hline
\multirow{2}{*}{$\kappa= \numprint{-1}$}  & $\mathcal{J}_{uc}$ & $\numprint{2.03e-2}$ & $\numprint{2.03e-2}$ & $\numprint{2.03e-2}$ & $\numprint{2.03e-2}$ & $\numprint{2.03e-2}$\\
 & $\mathcal{J}_c$ & $\numprint{1.93e-5}$ & $\numprint{8.17e-4}$ & $\numprint{1.45e-2}$ & $\numprint{2.02e-2}$ & $\numprint{2.03e-2}$\\
\hline
\end{tabular}
\end{small}
\caption{Source Control, No-Flux: Cost with $w=0$ ($\mathcal{J}_{uc}$) and with the optimal control ($\mathcal{J}_c$) for a range of $\kappa$ and $\beta$. Note that for $\beta = 10$, the cost functionals differ by $10^{-4}$, while for $\beta = 10^3$ they differ by $10^{-7}$.}
\label{TabSCN}
\end{table}

%% file: SCDExample2D.tex
\begin{table}
\centering
\begin{small}
\begin{tabular}{ | c | c || c | c | c | c | c |}
\hline
\multicolumn{2}{|c||}{}& $\beta = 10^{-5}$ & $\beta = 10^{-3}$ & $\beta = 10^{-1}$ & $\beta = 10^{1}$ & $\beta = 10^{3}$  \\
\hline
\hline
\multirow{2}{*}{$\kappa= \numprint{0}$}  & $\mathcal{J}_{uc}$ & $\numprint{1.50e-2}$ & $\numprint{1.50e-2}$ & $\numprint{1.50e-2}$ & $\numprint{1.50e-2}$ & $\numprint{1.50e-2}$\\
 & $\mathcal{J}_c$ & $\numprint{3.40e-5}$ & $\numprint{1.92e-3}$ & $\numprint{1.36e-2}$ & $\numprint{1.50e-2}$ & $\numprint{1.50e-2}$\\
\hline
\multirow{2}{*}{$\kappa= \numprint{1}$}  & $\mathcal{J}_{uc}$ & $\numprint{2.06e-2}$ & $\numprint{2.06e-2}$ & $\numprint{2.06e-2}$ & $\numprint{2.06e-2}$ & $\numprint{2.06e-2}$\\
 & $\mathcal{J}_c$ & $\numprint{4.27e-5}$ & $\numprint{2.49e-3}$ & $\numprint{1.85e-2}$ & $\numprint{2.06e-2}$ & $\numprint{2.06e-2}$\\
\hline
\multirow{2}{*}{$\kappa= \numprint{-1}$}  & $\mathcal{J}_{uc}$ & $\numprint{1.27e-2}$ & $\numprint{1.27e-2}$ & $\numprint{1.27e-2}$ & $\numprint{1.27e-2}$ & $\numprint{1.27e-2}$\\
 & $\mathcal{J}_c$ & $\numprint{2.88e-5}$ & $\numprint{1.61e-3}$ & $\numprint{1.18e-2}$ & $\numprint{1.27e-2}$ & $\numprint{1.27e-2}$\\
\hline
\end{tabular}
\end{small}
\caption{Source Control, Dirichlet:  Cost with $w=0$ ($\mathcal{J}_{uc}$) and with the optimal control ($\mathcal{J}_c$) for a range of $\kappa$ and $\beta$. Note that for $\beta = 10$ the cost functionals differ by $10^{-5}$, while for $\beta = 10^3$ they differ by $10^{-7}$ for $\kappa = 0$ and $\kappa = 1$, and by $10^{-8}$ for $\kappa = -1$.}
\label{TabSCD}
\end{table}

%% file: FCNExact2D.tex
\begin{table}
\resizebox{\textwidth}{!}{
\centering
\begin{small}
\begin{tabular}{ |c|| c || c | c | c | c | c |}
\hline
Solver& Error & $\beta = 10^{-5}$ & $\beta = 10^{-3}$ & $\beta = 10^{-1}$ & $\beta = 10^{1}$ & $\beta = 10^{3}$  \\
\hline
\hline
 Newton--Krylov& $\mathcal{E}_{\rho}$ & $\numprint{1.39e-15}$ & $\numprint{1.40e-14}$ & $\numprint{1.38e-13}$ & $\numprint{2.15e-13}$ & $\numprint{2.14e-13}$\\
 &$\mathcal{E}_{\adj}$ & $\numprint{1.82e-15}$ & $\numprint{1.82e-14}$ & $\numprint{1.81e-13}$ & $\numprint{4.44e-13}$ & $\numprint{4.45e-13}$\\
\hline
\hline
Fixed-Point& $\mathcal{E}_{\rho}$ & $\numprint{4.75e-8}$ & $\numprint{5.39e-8}$ & $\numprint{5.27e-8}$ & $\numprint{2.84e-8}$ & $\numprint{4.16e-8}$\\
& $\mathcal{E}_{\adj}$ & $\numprint{4.85e-8}$ & $\numprint{4.06e-8}$ & $\numprint{4.77e-8}$ & $\numprint{3.66e-8}$ & $\numprint{4.35e-8}$\\
\hline
\end{tabular}
\end{small}
}
\caption{Flow Control, No-Flux: Errors in computed solutions for state $\rho$ and adjoint $\adj$, compared with their analytic solutions, for a range of $\beta$.}
\label{TabFCNExact}
\end{table}